\documentclass[12pt]{article}
\usepackage{amsfonts}
\usepackage{amssymb}
\usepackage{amsmath}
=0pt
\def\sqr#1#2{{\vcenter{\vbox{\hrule height.#2pt
              \hbox{\vrule width.#2pt height#1pt \kern#1pt \vrule
width.#2pt}
              \hrule height.#2pt}}}}
\def\signed #1{{\unskip\nobreak\hfil\penalty50
              \hskip2em\hbox{}\nobreak\hfil#1
              \parfillskip=0pt \finalhyphendemerits=0 \par}}
\def\endpf{\signed {$\sqr69$}}
\def\e{\varepsilon}

\def\3n{\negthinspace \negthinspace \negthinspace }
\def\2n{\negthinspace \negthinspace }
\def\1n{\negthinspace }
\def\ns{\noalign{\smallskip} }

\def\ds{\displaystyle}
\def\D{\Delta}

\def\no{\noindent}

\def\ms{\medskip}
\def\bs{\bigskip}
\def\q{\quad}
\def\qq{\qquad}

\def\pa{\partial}

\def\wt{\widetilde}

\def\cds{\cdots}

\def\deq{\mathop{\buildrel\D\over=}}
\def\Re{{\mathop{\rm Re}\,}}

\def\({\Big (}
\def\){\Big )}
\def\[{\Big[}
\def\]{\Big]}

\def\={\buildrel \triangle \over =}

at 9pt 

\def\be{\begin{equation}}
\def\bel{\begin{equation}\label}
\def\ee{\end{equation}}
\def\bea{\begin{eqnarray}}
\def\eea{\end{eqnarray}}
\def\bt{\begin{theorem}}
\def\et{\end{theorem}}
\def\bc{\begin{corollary}}
\def\ec{\end{corollary}}
\def\bl{\begin{lemma}}
\def\el{\end{lemma}}
\def\bp{\begin{proposition}}
\def\ep{\end{proposition}}
\def\br{\begin{remark}}
\def\er{\end{remark}}
\def\ba{\begin{array}}
\def\ea{\end{array}}
\def\bd{\begin{definition}}
\def\ed{\end{definition}}

\newtheorem{lemma}{Lemma}[section]
\newtheorem{remark}{Remark}[section]

\newtheorem{theorem}{Theorem}[section]
\newtheorem{corollary}{Corollary}[section]

\newtheorem{definition}{Definition}[section]
\newtheorem{proposition}{Proposition}[section]

\makeatletter
   
   \@addtoreset{equation}{section}
\makeatother

\begin{document}

\title{\bf Carleman Estimates for Second Order Elliptic Operators with Limiting Weights, an Elementary Approach\thanks{This work is partially supported by the NSF of China
        under grants 12025105, 11971333 and 11931011, and by the Science Development Project of Sichuan University under
        grant 2020SCUNL201. }}

\author{Zengyu Li\thanks{School of Mathematics, Sichuan
University, Chengdu, 610064, China.  {\small\it
E-mail:} {\small\tt
lizengyu@stu.scu.edu.cn}.\ms} \q and \q Qi
L\"u\thanks{School of Mathematics, Sichuan
University, Chengdu, 610064, China.  {\small\it
E-mail:} {\small\tt lu@scu.edu.cn}.\ms} }

%\date{}
\maketitle

\begin{abstract}
By using some deep tools from microlocal
analysis, the authors of the papers (Ann. of
Math., 165 (2007), 567--591, J. Amer. Math.
Soc., 23 (2010), 655--691; Invent. Math., 178
(2009), 119--171; Duke Math. J., 158(2011),
83--120) have successfully established various
Carleman estimates for elliptic operators that
possess limiting Carleman weight. In this study,
we revisit these problems and present a unified
and fundamental approach for deriving these
estimates. The main tool we employ is an
elementary pointwise estimate for second-order
elliptic operators.
\end{abstract}

\no{\bf 2020 Mathematics Subject
Classification}. 35R30

\bs

\no{\bf Key Words}. Carleman estimate, limiting
Carleman weight, inverse problems.

%\newpage

\section{Introduction and main results}

In his groundbreaking paper \cite{C39}, T.
Carleman introduced a  revolutionary  method for
proving the strong unique continuation property
of second-order elliptic partial differential
equations (PDEs) in two variables. This method,
now known as the Carleman estimate, has since
become a fundamental tool in the study of
various important problems in PDEs, including
unique continuation problems, inverse problems,
and control problems.

In recent years, the Carleman estimate has been
successfully applied to solve the famous
Calder\'{o}n problem, and several deep Carleman
estimates with limiting weight functions have
been established
\cite{DKSU09,HT13,GT11,IUY10,KSU07}. The proofs
of these Carleman estimates in
\cite{DKSU09,HT13,GT11,KSU07} relied on
sophisticated techniques from microlocal
analysis, such as the Fefferman-Phong
inequality. Additionally, separate treatments
were required to handle scenarios where the
space dimension is either $2$ or greater than or
equal to $3$.

In this paper, we aim to provide a unified and
elementary approach to derive the Carleman
estimates in
\cite{DKSU09,HT13,GT11,IUY10,KSU07}. With this
approach, we are able to provide a unified
treatment for both $2$-dimensional and
higher-dimensional cases. We believe that our
unified approach not only simplifies the proofs
of these Carleman estimates but also provides
some new insights into the underlying theory.

To present the main results of this paper, we
begin by revisiting the notations and
definitions used for Riemannian manifolds. More
comprehensive explanations can be found in the
reference \cite{J17}.

Let $M$ be an $n$-dimensional $C^{3}$-smooth
compact Riemannian manifold  with a
$C^{2}$-smooth boundary. In this context, we
will adopt the following notations: $g$
represents the $C^3$-smooth Riemannian metric
tensor on $M$, $\langle\cdot,\cdot\rangle_g$ and
$|\cdot|_g$ denote the inner product and norm on
the tangent vector fields with respect to $g$
respectively. The Levi-Civita connection induced
by $g$ on $M$ is denoted by $\mathcal{D}_{g}$.
The gradient operator, divergence operator,
Hesse operator, and Laplace-Beltrami operator on
$M$ will be denoted by $\nabla_g$,
$\mathrm{div}_g$, $\mathrm{Hess}_g$   and
$\Delta_g$, respectively. $dV_g$? denotes the
volume form on $(M,g)$, while $dS_g$ signifies
the induced volume form on $\partial M$. When
working with the Euclidean metric
$\mathfrak{e}$, we will omit the subscripts of
the inner product, the norm, the operators for
simplicity. Denote by $dx$? the volume form on
$(M,\mathfrak{e})$, while $dS$ signifies the
induced volume form on $\partial M$.

% The $L^{2}$ norm of a function is given by
%\begin{equation*}
%   \|u\|_{L^{2}(M)}=\Big(\int_{M}|u|^{2}\,dV_{g}\Big)^{\frac{1}{2}},
%\end{equation*}
%and the norm on boundary is given by
%\begin{equation*}
%   \|u\|_{L^{2}(\partial M)}=\Big(\int_{\partial M}|u|^{2}\,dS_{g}\Big)^{\frac{1}{2}}.
%\end{equation*}
%We write for short
%\begin{equation*} \|\nabla_{g}u\|_{L^{2}(M)}=\||\nabla_{g}u|_{g}\|_{L^{2}(M)}=\Big(\int_{M}|\nabla_{g}u|_{g}^{2}\,dV_{g}\Big)^{\frac{1}{2}}.
%\end{equation*}

Let $(N,g)$ be a $n$-dimensional  $C^{3}$-smooth
open  Riemannian manifold  such that
$M\subset\subset N$. Let us recall the
definition of the limiting Carleman weight for
the Laplace-Beltrami operator $\Delta_g$ on $N$.

\begin{definition}\label{defn1.1}
Let $\varphi\in C^{3}(N;\mathbb{R})$, it is
called a limiting Carleman weight for the
Laplace-Beltrami operator if it has
non-vanishing gradient, and satisfies
\begin{equation}\label{1.1} \mathrm{Hess}_{g}\,\varphi(X,X)+\mathrm{Hess}_{g}\,\varphi(\nabla_{g}\varphi,\nabla_{g}\varphi)=0 \q \mbox{ in }\ N
\end{equation}
for all $X\in T(N)$ satisfying
$|X|_{g}^{2}=|\nabla_{g}\varphi|_{g}^{2}$ and
$\langle X,\nabla_{g}\varphi\rangle_{g}=0$.
\end{definition}
\begin{remark}
The initial definition of limiting Carleman
weights, as originally stated in \cite{KSU07},
utilizes terminology from semiclassical
analysis.  In our current paper, we adopt an
alternative definition that relies solely on
elementary concepts. The proof of the
equivalence between these two definitions can be
found in \cite{DKSU09}.
\end{remark}
\begin{remark}
As mentioned in \cite{DKSU09}, it is well-known
that a generic manifold in dimension $n\geq3$
may not possess a limiting Carleman weight.
Consequently, our focus in this paper is on
Riemannian manifolds that do have such weights.
In \cite{AFG17}, certain conditions for the
existence of limiting Carleman weights on a
manifold have been investigated, which are
closely related to the properties of the Weyl
tensor and the Cotton-York tensor. In order for
these tensors to exist, it is assumed that the
metric tensor $g$ is at least $C^{3}$-smooth.

However, the situation is different in the case
of dimension $n=2$. Here, any harmonic function
with a non-vanishing gradient can be considered
as a limiting Carleman weight.
\end{remark}
\begin{remark}
All manifolds throughout this paper are assumed
to be oriented and connected.
\end{remark}

We have the following two Carleman estimates,
both of which are proved by a fundamental
pointwise identity in Section 2. The first one
is a Carleman estimate with limiting weights on
a Riemannian manifold of dimension $n\geq3$.

\begin{theorem}\label{thm1.1}
Assume that $n\geq3$. Let $(N,g)$ be an
$n$-dimensional $C^{3}$-smooth open Riemannian
manifold, and $(M,g)$ an $n$-dimensional
$C^{3}$-smooth compact manifold  with a
$C^{2}$-smooth boundary, such that
$M\subset\subset N$. Suppose that $\varphi$ is a
limiting Carleman weight on $(N,g)$. Let $X$ be
an $L^{\infty}$ vector field on $M$ and $q\in
L^{\infty}(M)$. Let $\nu$ be the outward unit
normal vector field along $\partial M$. Then
there exist  two constants $C>0$ and
$\tau_{0}>0$ such that for $\tau\geq\tau_{0}$
and for  all functions $v\in H^{2}(N)$, we have
\begin{equation}\label{1.2}
\begin{aligned}
&\tau^{3}\big\|e^{\tau\varphi}v\big\|_{L^{2}(\partial
M)}^{2}
+\tau\big\|e^{\tau\varphi}\nabla_{\parallel}v\big\|_{L^{2}(\partial
M)}^{2}
+\tau^{2}\int_{\partial M}e^{2\tau\varphi}|\nabla_{\perp}v|_{g}|v|\,dS_{g}\\
&+\tau\int_{\partial
M}e^{2\tau\varphi}|\nabla_{\parallel}v|_{g}|\nabla_{\perp}v|_{g}\,dS_{g}
+\tau\int_{\partial M}e^{2\tau\varphi}|\langle\nabla_{\perp}\varphi,\nu\rangle_{g}|\,|\nabla_{\perp} v|_{g}^{2}\,dS_{g}\\
&+\big\|e^{\tau\varphi}(-\Delta_{g}+X+q)v\big\|_{L^{2}(M)}^{2}
\geq
C(\tau^{2}\big\|e^{\tau\varphi}v\big\|_{L^{2}(M)}^{2}+\big\|e^{\tau\varphi}\nabla_{g}
v\big\|_{L^{2}(M)}^{2}),
\end{aligned}
\end{equation}
where
$\nabla_{\perp}v=\langle\nabla_{g}v,\nu\rangle_{g}\nu$,
and
$\nabla_{\parallel}v=\nabla_{g}v-\nabla_{\perp}v$.
\end{theorem}

%\begin{remark}\label{rem1.4}
%Theorem \ref{thm1.1} is a Carleman estimate on Riemannian manifolds. When  $\widetilde{\Omega}\subset\mathbb{R}^{n}$ is a bounded domain, and $\Omega\subset\subset\widetilde{\Omega}$ is an open set with a $C^{2}$ boundary, then  \eqref{1.2} is already proved in \cite{KSU07} by employing semiclassical analysis. We can also directly get a Carleman estimate for $v\in H^{2}(\Omega)$ from \eqref{1.2} which is first used in \cite{KkU16} in the case of conductivities with lower regularity, and the potential may be possibly $\tau$-dependent. We will mention these two Carleman estimates as corollaries of Theorem \ref{thm1.1} in Section 3. The main difference between these two estimates lies in the estimation of boundary terms, and the latter requires the expression for the Laplacian and the gradient on the boundary. As we will see later, less computation is required in our proof since we can get rid of the second order derivatives in boundary terms by the fundamental pointwise identity.
%\end{remark}

Here and in what follows, we denote by $C$ a
generic constant, which may vary from line to
line. When we want to distinguish several
constants, we use the notaions $C_1$,
$C_2,\cds$, etc.

%
%\begin{remark}
%When the manifold, the boundary, and the metric tensor are all of $C^{\infty}$-smooth, we can derive from \eqref{1.2} a Carleman estimate which is first established in \cite{DKSU09} to deal with anisotropic Calder\'{o}n inverse problems.  The original proof in \cite{DKSU09} is also based on symbol calculus.
%\end{remark}
%

From Theorem \ref{thm1.1},  we can obtain the
following four Carleman estimates. In the
following we replace $\tau$ by $\frac{1}{h}$ for
small $h$ in order to follow the notations in
\cite{DKSU09,KS13,KSU07,KkU16}.
\begin{corollary}\cite[Theorem 4.1]{DKSU09}\label{cor1.1}
Under the same conditions of Theorem
\ref{thm1.1}, there exist two constants
$C,h_{0}>0$ such that for $0<h\leq h_{0}$, one
has
\begin{equation}\label{1.3}
\big\|e^{\frac{\varphi}{h}}v\big\|_{L^{2}(M)}^{2}
+\big\|e^{\frac{\varphi}{h}}h\nabla_{g}v\big\|_{L^{2}(M)}^{2}\leq
Ch^{2}\big\|e^{\frac{\varphi}{h}}(-\Delta_{g}+X+q)v\big\|_{L^{2}(M)}^{2}
\end{equation}
for any $v\in H_{0}^{2}(M)$.
\end{corollary}
\begin{remark}
If $(N,g)$ is a $C^\infty$-smooth open
Riemannian manifold, $(M,g)$ a $C^\infty$-smooth
compact Riemannian submanifold with boundary
such that $M\subset\subset N$,   $X$ is a smooth
vector field on $M$ and $q$ is a
$C^\infty$-smooth function on $M$,  the
inequality \eqref{1.3} is proved in
\cite{DKSU09} for $v\in C_{0}^{\infty}(M)$ (see
\cite[Theorem 4.1]{DKSU09}).
\end{remark}

To present the next result, we first  recall the
concept of admissible manifold.
\begin{definition}\label{defn1.2}
Let $(M,g)$ be a $n$-dimensional $C^{3}$-smooth
compact Riemannian manifold with the
$C^{2}$-smooth boundary, and $n\geq3$. We say
that $(M,g)$ is admissible if it satisfies

(1) $(M,g)\subset\subset (\mathbb{R}\times
M_{0},g)$, and $g=c(\mathfrak{e}\oplus g_{0})$,
where $(M_{0},g_{0})$ is a compact
$(n-1)$-dimensional manifold with $C^{2}$-smooth
boundary, $\mathfrak{e}$ is the Euclidean metric
on the real line, and $c$ is a smooth positive
function in the cylinder $\mathbb{R}\times
M_{0}$.

(2) $(M_{0},g_{0})$ is simple, i.e., $\partial
M_{0}$ is strictly convex and for any $p\in
M_{0}$, the exponential map $\exp_{p}$ is a
diffeomorphism from its maximal domain of
definition in $T_{p}M_{0}$ onto $M_{0}$.
\end{definition}
\begin{remark}
The motivation to introduce the admissible
manifold  lies in that there exist  limiting
Carleman weights on such kind of manifold (see
\cite{DKSU09}). Classical examples for
admissible manifold include bounded domains in
Euclidean space, in the sphere minus a point,
and in Hyperbolic space.
\end{remark}

If $(M,g)$ is admissible, then points of $x\in
M$ can be written as $x=(x_{1},x')$, where
$x_{1}$ is the Euclidean coordinate. We define
\begin{equation*}
\begin{cases}\ds
\partial M_{\pm} =\{x\in\partial M:\pm\partial_{\nu}\varphi(x)>0\},\\
\ns\ds \partial M_{\tan} =\{x\in\partial
M:\partial_{\nu}\varphi(x)=0\},
\end{cases}
\end{equation*}
where $\varphi(x)=x_{1}$ is a natural limiting
Carleman weight on $(M,g)$.

\begin{corollary}\label{cor1.2}
Let $(M,g)$ be admissible, $q\in L^{\infty}(M)$
and $\varphi(x)=\pm x_{1}$. Denote by
$\partial_{\nu}$ the outward unit normal vector
field to $\partial M$. Then there exist two
constants $C,h_{0}>0$ such that for $0<h\leq
h_{0}$ and $\delta>0$, one has
\begin{equation}\label{1.4}
\begin{aligned}
&\delta
h^{3}\big\|\partial_{\nu}u\big\|_{L^{2}(\{\partial_{\nu}\varphi\leq-\delta\})}^{2}
+h^{4}\big\|\partial_{\nu}u\big\|_{L^{2}(\{-\delta<\partial_{\nu}\varphi<h/3\})}^{2}
+h^{2}\big(\big\|u\big\|_{L^{2}(M)}^{2}+\big\|h\nabla_{g}u\big\|_{L^{2}(M)}^{2}\big) \\
&\leq
C\big(\big\|e^{\frac{\varphi}{h}}(-h^{2}\Delta_{g}+h^{2}q)(e^{-\frac{\varphi}{h}}u)\big\|_{L^{2}(M)}^{2}
+h^{3}\big\|\partial_{\nu}u\big\|_{L^{2}(\{\partial_{\nu}\varphi\geq
h/3\})}^{2} \big)
\end{aligned}
\end{equation}
for any $u\in H^2(M)\cap H^1_0(M)$.
\end{corollary}
\begin{remark}
If $(M,g)$ is a $C^\infty$-smooth compact
Riemannian  manifold the inquality \eqref{1.4}
is proved in \cite{KS13} for $v\in
C^{\infty}(M)$ (see \cite[Proposition
4.2]{KS13}).
\end{remark}
\begin{corollary}\label{cor1.3}
Let $\widetilde{\Omega}\subset\mathbb{R}^{n}$,
$n\geq 3$ be a bounded open set. Let
$\Omega\subset\subset\widetilde{\Omega}$ be an
open set with a $C^{2}$-smooth boundary
$\partial\Omega$. Suppose that $\varphi\in
C^{3}(\widetilde{\Omega})$ is a limiting
Carleman weight. Let $q\in L^{\infty}(\Omega)$.
Denote by $\nu$ the unit outward normal vector
to $\partial\Omega$ and define
\begin{equation*}
\partial \Omega_{\pm}=\big\{x\in\partial\Omega:\pm\partial_{\nu}\varphi(x)\geq 0\big\}.
\end{equation*}
Then there exist two constants $C,h_{0}>0$ such
that for $0<h\leq h_{0}$, one has
\begin{equation}\label{1.5}
\begin{aligned}
&-\frac{h^3}{C}\int_{\partial\Omega_{-}}\partial_{\nu}\varphi\big|e^{\frac{\varphi}{h}}\partial_{\nu}v\big|^{2}\,dS
+\frac{h^{2}}{C}\big(\big\|e^{\frac{\varphi}{h}}v\big\|^{2}_{L^{2}(\Omega)}
+\big\|e^{\frac{\varphi}{h}}h\nabla v\big\|^{2}_{L^{2}(\Omega)}\big)\\
&\leq
\big\|e^{\frac{\varphi}{h}}(-h^{2}\Delta+h^{2}q)v\big\|^{2}_{L^{2}(\Omega)}
+Ch^{3}\int_{\partial\Omega_{+}}\partial_{\nu}\varphi\big|e^{\frac{\varphi}{h}}\partial_{\nu}v\big|^{2}\,dS
\end{aligned}
\end{equation}
for any $v\in H^2(\Omega)\cap H^1_0(\Omega)$.
\end{corollary}
\begin{remark}
If $\pa\Omega$ is  $C^\infty$-smooth and
$\varphi\in C^{\infty}(\widetilde{\Omega})$, the
inquality \eqref{1.5} is proved in \cite{KSU07}
for $v\in C^{\infty}(M)$ (see \cite[Proposition
3.2]{KSU07}).
\end{remark}

\begin{corollary}\label{cor1.4}
Let $\widetilde{\Omega}\subset\mathbb{R}^{n}$,
$n\geq 3$ be a bounded open set. Let
$\Omega\subset\subset\widetilde{\Omega}$ be an
open set with $C^{2}$-smooth boundary. Consider
the operator
\begin{equation*}
-\Delta+\mathcal{A}\cdot\nabla+q
\end{equation*}
where $\mathcal{A}\in
L^{\infty}(\Omega;\mathbb{C}^{n})$, $q\in
L^{\infty}(\Omega;\mathbb{C})$ are possibly
$h$-dependent with
\begin{equation*}
\|\mathcal{A}\|_{L^{\infty}(\Omega)}=\mathcal{O}(1),
\quad
\|q\|_{L^{\infty}(\Omega)}=\mathcal{O}\(\frac{1}{h}\)
\end{equation*}
as $h\rightarrow 0$. Suppose that $\varphi\in
C^{3}(\widetilde{\Omega})$ is a limiting
Carleman weight. Let $\nu$ denote the unit
outward normal vector to $\partial\Omega$, and
$\partial\Omega_{\pm}$ be as in Corollary
\ref{cor1.3}. Then there exist two constants
$C,h_{0}>0$ such that for $0<h\leq h_{0}$, one
has
\begin{equation}\label{1.6}
\begin{aligned}
&h\big\|e^{\frac{\varphi}{h}}v\|_{L^{2}(\partial\Omega)}^{2}
+h^{2}\int_{\partial\Omega}e^{\frac{2\varphi}{h}}|\partial_{\nu}v||v|\,dS
+h^{3}\big\|e^{\frac{\varphi}{h}}\nabla_{t}v\big\|_{L^{2}(\partial\Omega)}^{2}
+h^{3}\int_{\partial\Omega}e^{\frac{2\varphi}{h}}|\nabla_{t}v||\partial_{\nu}v|\,dS\\
&-h^{3}\int_{\partial\Omega_{-}}\partial_{\nu}\varphi\big|e^{\frac{\varphi}{h}}\partial_{\nu} v\big|^{2}\,dS+\big\|e^{\frac{\varphi}{h}}(-h^{2}\Delta+h\mathcal{A}\cdot h\nabla+h^{2}q)v\big\|_{L^{2}(\Omega)}^{2}\\
\geq&C\Big[
h^{2}\big(\big\|e^{\frac{\varphi}{h}}v\big\|_{L^{2}(\Omega)}^{2}+\big\|e^{\frac{\varphi}{h}}h\nabla
v\big\|_{L^{2}(\Omega)}^{2}\big)
+h^{3}\int_{\partial\Omega_{+}}\partial_{\nu}\varphi\big|e^{\frac{\varphi}{h}}\partial_{\nu}v\big|^{2}\,dS\Big]
\end{aligned}
\end{equation}
for any $v\in H^{2}(\Omega)$, where $\nabla_{t}$
denote the tangential component of the gradient.
\end{corollary}
\begin{remark}
If $\varphi\in C^{\infty}(\widetilde{\Omega})$,
the inquality \eqref{1.6} is proved in
\cite{KkU16} (see \cite[Proposition
3.2]{KkU16}).
\end{remark}

In prior literature,  Corollary \ref{cor1.1} has
been employed in solving anisotropic
Calder\'{o}n problems. Corollaries \ref{cor1.2}
and  \ref{cor1.3} are applied to investigate
Calder\'{o}n problems with partial data in
dimensions $n\geq3$, while Corollary
\ref{cor1.4} is applicable to situations
involving less regular conductivities.

The Carleman estimate for the case $n=2$ is
established specifically for Riemann surfaces.
We restrict our attention to the case when the
surface is simply connected. In this scenario,
we select a weight function that is a harmonic
Morse function. However, this particular choice
of weight function may introduce some critical
points, which renders it no longer a limiting
weight. Instead, it is referred to as a
degenerate weight.

\begin{theorem}\label{thm1.2}
Let $(\widetilde{N},\widetilde{g})$ be a compact
connected Riemann surface, and
$(\widetilde{M},\widetilde{g})$ be a compact
connected Riemann surface with boundary such
that $\widetilde{M}\subset\widetilde{N}$, where
$\widetilde{g}$ is the $C^{\infty}$-smooth
metric tensor. Let
$\varphi:\widetilde{N}\rightarrow\mathbb{R}$ be
a harmonic Morse function with prescribed
critical points $\{p_{1},p_{2},\cdots,p_{m}\}$
in the interior of $\widetilde{M}$, and critical
points $\{q_{1},\cdots,q_{s}\}$ on
$\partial\widetilde{M}$. Denote by
$\partial_{\nu}$ the outward unit normal vector
field to $\partial \widetilde{M}$. Define
$\Gamma_{0}=\{x\in\partial
\widetilde{M}:\partial_{\nu}\varphi(x)=0\}$, and
let $\Gamma=\partial
\widetilde{M}\setminus\Gamma_{0}$ be its
complement. Then for all $q\in
L^{\infty}(\widetilde{M})$, there exists two
constants $C>0$ and $\tau_{0}>0$ such that for
all functions $v\in H^{2}(\widetilde{M})\cap
H_{0}^{1}(\widetilde{M})$, we have for
$\tau\geq\tau_{0}$,
\begin{equation}\label{1.7}
\begin{aligned} &\tau\big\|e^{\tau\varphi}v\big\|_{L^{2}(\widetilde{M})}^{2}+\big\|e^{\tau\varphi}v\big\|_{H^{1}(\widetilde{M})}^{2}
+\tau^{2}\big\||\nabla_{\widetilde{g}}\varphi|_{\widetilde{g}}e^{\tau\varphi}v\big\|_{L^{2}(\widetilde{M})}^{2}+
\big\|e^{\tau\varphi}\partial_{\nu}v\big\|_{L^{2}(\Gamma_{0})}^{2}\\
\leq
&C\Big(\big\|e^{\tau\varphi}(-\Delta_{\widetilde{g}}+q)v\big\|_{L^{2}(\widetilde{M})}^{2}
+\tau\big\|e^{\tau\varphi}\partial_{\nu}v\big\|_{L^{2}(\Gamma)}^{2}\Big).
\end{aligned}
\end{equation}
\end{theorem}
%
%\begin{remark}\label{rem1.6}
%This result has an early version in \cite{GT11}. As a comparison, the $C^{\infty}$-regularity of the boundary in \cite{GT11} is weakened to be $C^{2}$ in Theorem \ref{1.2}.
%\end{remark}
%

In the final part of this section, as a direct
result of Theorem \ref{thm1.2}, we provide a
Carleman estimate that is utilized in the
solution of the two-dimensional Calder\'{o}n
problem with partial data.

To begin with, we introduce some notations from
\cite{IUY10}. Let $\Omega\subset\mathbb{R}^{2}$
be a bounded domain with a smooth boundary. Let
$\Gamma_{1}\subset\partial\Omega$ be a nonempty
open subset of the boundary, and
$\Gamma_{2}=\partial\Omega\setminus\overline{\Gamma_{1}}$.
We identify $x=(x_{1},x_{2})\in\mathbb{R}^{2}$
with $z=x_{1}+ix_{2}\in\mathbb{C}$. We use the
notations
$\partial_{z}=\frac{1}{2}(\partial_{x_{1}}-i\partial_{x_{2}})$,
and
$\partial_{\overline{z}}=\frac{1}{2}(\partial_{x_{1}}+i\partial_{x_{2}})$.
Let
$\Phi(z)=\varphi(x_{1},x_{2})+i\psi(x_{1},x_{2})\in
C^{2}(\overline{\Omega})$ be holomorphic in
$\Omega$, that is
\begin{equation}\label{1.8}
\partial_{\overline{z}}\Phi(z)=0\quad\text{ in }\,\Omega.
\end{equation}
Denote by $\mathcal{H}$ the set of critical
points of $\Phi$, that is,
\begin{equation*}
\mathcal{H}=\{z\in\overline{\Omega}:\partial_{z}\Phi(z)=0\}.
\end{equation*}
Assume that $\Phi$ has no critical points on
$\overline{\Gamma_{1}}$, and all the critical
points are non-degenerate, i.e.,
\begin{equation}\label{1.9}
\partial_{z}^{2}\Phi(z)\neq0,\quad \forall z\in\mathcal{H}.
\end{equation}
We also assume that $\Phi$ satisfies
\begin{equation}\label{1.10}
\Gamma_{2}\subset\{x\in\partial\Omega:\partial_{\nu}\varphi(x)=0\}.
\end{equation}

It follows immediately from Theorem \ref{thm1.2}
that
\begin{corollary}\cite[Proposition 5.3]{IUY10}\label{cor1.5}
Suppose that $\Phi$ satisfies
\eqref{1.8}-\eqref{1.10}. Let $v\in
H^{2}(\Omega)\cap H_{0}^{1}(\Omega)$ be a
real-valued function. Denote by $\nu$ the unit
outward normal vector to $\partial\Omega$. Then
there exist  two constants $C>0$ and
$\tau_{0}>0$ such that for all $|\tau|\geq
\tau_{0}$, we have
\begin{equation}\label{1.11}
\begin{aligned}
&|\tau|
\big\|e^{\tau\varphi}v\big\|^{2}_{L^{2}(\Omega)}
+\big\|e^{\tau\varphi}v\big\|^2_{H^{1}(\Omega)}
+\big\|e^{\tau\varphi}\partial_{\nu}v\big\|^{2}_{L^{2}(\Gamma_{2})}
+\tau^{2}\big\||\partial_{z}\Phi|e^{\tau\varphi}v\big\|^{2}_{L^{2}(\Omega)}\\
&\leq  C\Big(\big\|e^{\tau\varphi}\Delta
v\big\|^{2}_{L^{2}(\Omega)}+|\tau|\int_{\Gamma_{1}}|e^{\tau\varphi}\partial_{\nu}v|^{2}\,dS\Big).
\end{aligned}
\end{equation}
\end{corollary}
%

%\begin{remark}\label{rem1.7}
%Theorem \ref{thm1.2} is a general result on Riemann surfaces. Now we consider the special case of Euclidean spaces, i.e.,   $\widetilde{\Omega}\subset\mathbb{R}^{2}$ is a bounded open set, and $\Omega\subset\subset\widetilde{\Omega}$ is an open set with $C^{2}$ boundary. let $\Omega\subset\mathbb{R}^{2}$ be a bounded open set. We further consider the special case that $\Omega$ has a $C^{\infty}$-smooth boundary, then we can obtain a Carleman estimate for $v\in H^{2}(\Omega)\cap H_{0}^{1}(\Omega)$ from \eqref{1.3}. This Carleman estimate firstly appeared in \cite{IUY10} is used in the Calder\'{o}n problem with partial data in two dimensions. The original proofs in \cite{GT11,IUY10} only require some elementary calculations, especially the integration by parts. Anyhow, our aim in this paper is to develop a unified approach to derive the Carleman estimates for both $n=2$ and $n\geq3$, and our proof is slightly different.
%\end{remark}

%The Carleman estimates with limiting weights and degenerate weights have two main applications in the Calder\'{o}n problem with partial data. One is to give a new method for constructing the complex geometrical optics (CGO) solutions, the other is to prove some uniqueness results. Our results of Carleman estimates weaken the condition of the boundary regularity to be of $C^{2}$ smooth.

The remainder of this paper is structured as
follows.

In Section 2, we present a pivotal pointwise
weighted identity for the Laplace-Beltrami
operators on Riemannian manifolds, which forms
the foundation for deriving the aforementioned
Carleman estimates. With the aid of this
identity, we establish Theorems \ref{thm1.1} and
\ref{thm1.2} in Sections 3 and 4, respectively.

In Section 5, we delve into the discussion of
some results concerning the Calder\'{o}n problem
with partial data through the utilization of
Carleman estimates.

\section{A fundamental weighted identity}

In this section, we establish a fundamental
weighted identity for Laplace-Beltrami operators
on Riemannian manifolds.

Let $(\mathcal{M},\mathfrak{g})$ be a
$C^{3}$-smooth Riemannian manifold of dimension
$n$ with a $C^{3}$-smooth metric tensor
$\mathfrak{g}$. The meaning of
$\langle\cdot,\cdot\rangle_{\mathfrak{g}}$,
$\mathcal{D}_{\mathfrak{g}}$,
$\nabla_{\mathfrak{g}}$,
$\mathrm{div}_{\mathfrak{g}}$,
$\mathrm{Hess}_{\mathfrak{g}}$,
$\Delta_{\mathfrak{g}}$ can be understood as
mentioned in Section 1. Let $X,Y$ be
$C^{1}$-smooth vector fields on $\mathcal{M}$
and $f\in C^{1}(\mathcal{M})$. We first recall
the following results.
\begin{eqnarray}
\mathrm{div}_{\mathfrak{g}}\,(fX)&=&\langle \nabla_{\mathfrak{g}}f,X\rangle_{\mathfrak{g}}+f\mathrm{div}_{\mathfrak{g}}\,X, \label{2.1}\\
\nabla_{\mathfrak{g}}\langle
X,Y\rangle_{\mathfrak{g}}&=&(\mathcal{D}_{\mathfrak{g}}X,Y)_{\mathfrak{g}}+(X,\mathcal{D}_{\mathfrak{g}}Y)_{\mathfrak{g}},
\label{2.2}
\end{eqnarray}
where
$(\mathcal{D}_{\mathfrak{g}}X,Y)_{\mathfrak{g}}$
stands for the contraction of $\mathfrak{g}
\otimes \mathcal{D}_{\mathfrak{g}}X \otimes Y$.

For $v\in C^{2}(\mathcal{M})$, fix $\ell\in
C^{3}(\mathcal{M})$, put
\begin{equation*}
\theta=e^{\ell},\quad u=\theta v.
\end{equation*}
From \eqref{2.1} we have
$$
\begin{aligned}
\theta\Delta_{\mathfrak{g}}v&=\theta\mathrm{div}_{\mathfrak{g}}\,[\nabla_{\mathfrak{g}}(\theta^{-1}u)]\\
&=\theta\mathrm{div}_{\mathfrak{g}}\,(-\theta^{-1}u\nabla_{\mathfrak{g}}\ell+\theta^{-1}\nabla_{\mathfrak{g}}u)\\
&=\Delta_{\mathfrak{g}}u-2\langle\nabla_{\mathfrak{g}}\ell,\nabla_{\mathfrak{g}}u\rangle_{\mathfrak{g}}
+|\nabla_{\mathfrak{g}}\ell|_{\mathfrak{g}}^{2}u-(\Delta_{\mathfrak{g}}\ell)u.
\end{aligned}
$$
Choose a symmetric matrix
$Q=Q(x)=(q_{j}^{i}(x))_{n\times n}$ with
$q_{j}^{i}\in C^{1}(\mathcal{M})$. Denote by $I$
the unit matrix of order $n$. Put
$\theta\Delta_{\mathfrak{g}}v=I_{1}+I_{2}$,
where
\begin{equation}\label{2.3}\left\{
\begin{aligned}
I_{1}&=\Delta_{\mathfrak{g}}u+\mathrm{div}_{\mathfrak{g}}\,(Q\nabla_{\mathfrak{g}}u)+\langle Y,\nabla_{\mathfrak{g}}u\rangle_{\mathfrak{g}}+(|\nabla_{\mathfrak{g}}\ell|_{\mathfrak{g}}^{2}+R)u\\
&=\mathrm{div}_{\mathfrak{g}}\,\big[(Q+I)\nabla_{\mathfrak{g}}u\big]+\langle Y,\nabla_{\mathfrak{g}}u\rangle_{\mathfrak{g}}+(|\nabla_{\mathfrak{g}}\ell|_{\mathfrak{g}}^{2}+R)u,\\
I_{2}&=-\mathrm{div}_{\mathfrak{g}}\,(Q\nabla_{\mathfrak{g}}u)-\langle 2\nabla_{\mathfrak{g}}\ell+Y,\nabla_{\mathfrak{g}}u\rangle_{\mathfrak{g}}-(\Delta_{\mathfrak{g}}\ell+R)u,\\
\end{aligned}
\right.
\end{equation}
where $Y$ is a $C^{1}$-smooth vector field on
$\mathcal{M}$ and $R\in C^{1}(\mathcal{M})$.

We have the following pointwise identity, which
is in fact true over general semi-Riemannian
manifolds.
\begin{theorem}\label{thm2.1}
It holds that
\begin{equation}\label{2.4}
|e^{\ell}\Delta_{\mathfrak{g}}v|^{2}+\mathrm{div}_{\mathfrak{g}}\,V=|I_{1}|^{2}+|I_{2}|^{2}+B_{1}u^{2}+2B_{2}u+2(B_{3}+F)+B_{4},
\end{equation}
where
\begin{equation}\label{2.5}
\left\{
\begin{aligned}
B_{1}=&(\Delta_{\mathfrak{g}}
\ell+R)\mathrm{div}_{\mathfrak{g}}Y+(|\nabla_{\mathfrak{g}}\ell|_{\mathfrak{g}}^{2}+R)\mathrm{div}_{\mathfrak{g}}(2\nabla_{\mathfrak{g}}\ell+Y)
+\langle\nabla_{\mathfrak{g}}(\Delta_{\mathfrak{g}}\ell+R),Y \rangle_{\mathfrak{g}}\\
&+\langle\nabla_{\mathfrak{g}}(|\nabla_{\mathfrak{g}}\ell|_{\mathfrak{g}}^{2}+R),2\nabla_{\mathfrak{g}}\ell+Y
\rangle_{\mathfrak{g}}
-2(|\nabla_{\mathfrak{g}}\ell|_{\mathfrak{g}}^{2}+R)(\Delta_{\mathfrak{g}}\ell+R),\\
B_{2}=&\langle\nabla_{\mathfrak{g}}(\Delta_{\mathfrak{g}}\ell+R),(Q+I)\nabla_{\mathfrak{g}}u\rangle_{\mathfrak{g}}
+\langle\nabla_{\mathfrak{g}}(|\nabla_{\mathfrak{g}}\ell|_{\mathfrak{g}}^{2}+R),Q\nabla_{\mathfrak{g}}u\rangle_{\mathfrak{g}},\\
B_{3}=&\langle(\mathcal{D}_{\mathfrak{g}}(2\nabla_{\mathfrak{g}}\ell+Y),\nabla_{\mathfrak{g}}u)_{\mathfrak{g}},(Q+I)\nabla_{\mathfrak{g}}u\rangle_{\mathfrak{g}}
+\langle(\mathcal{D}_{\mathfrak{g}}Y,\nabla_{\mathfrak{g}}u)_{\mathfrak{g}},Q\nabla_{\mathfrak{g}}u\rangle_{\mathfrak{g}}\\
&+(\Delta_{\mathfrak{g}}\ell+R)\langle\nabla_{\mathfrak{g}}u,(Q+I)\nabla_{\mathfrak{g}}u\rangle_{\mathfrak{g}}
+(|\nabla_{\mathfrak{g}}\ell|^{2}_{\mathfrak{g}}+R)\langle\nabla_{\mathfrak{g}}u,Q\nabla_{\mathfrak{g}}u\rangle_{\mathfrak{g}}\\
&-\langle Y,\nabla_{\mathfrak{g}}u\rangle_{\mathfrak{g}}\langle2\nabla_{\mathfrak{g}}\ell+Y,\nabla_{\mathfrak{g}}u\rangle_{\mathfrak{g}},\\
B_{4}=&-2\mathrm{div}_{\mathfrak{g}}\big[(Q+I)\nabla_{\mathfrak{g}}u\big]\mathrm{div}_{\mathfrak{g}}(Q\nabla_{\mathfrak{g}}u),\\
F=&\langle(2\nabla_{\mathfrak{g}}\ell+Y,\mathcal{D}_{\mathfrak{g}}\nabla_{\mathfrak{g}}u)_{\mathfrak{g}},(Q+I)\nabla_{\mathfrak{g}}u\rangle_{\mathfrak{g}}
+\langle(Y,\mathcal{D}_{\mathfrak{g}}\nabla_{\mathfrak{g}}u)_{\mathfrak{g}},Q\nabla_{\mathfrak{g}}u\rangle_{\mathfrak{g}},\\
V=&2\big[\langle2\nabla_{\mathfrak{g}}\ell+Y,\nabla_{\mathfrak{g}}u\rangle_{\mathfrak{g}}+(\Delta_{\mathfrak{g}}\ell+R)u\big](Q+I)\nabla_{\mathfrak{g}}u
\\
&+2\big[\langle Y,\nabla_{\mathfrak{g}}u\rangle_{\mathfrak{g}}+(|\nabla_{\mathfrak{g}}\ell|_{\mathfrak{g}}^{2}+R)u\big]Q\nabla_{\mathfrak{g}}u\\
&+(|\nabla_{\mathfrak{g}}\ell|_{\mathfrak{g}}^{2}+R)u^{2}(2\nabla_{\mathfrak{g}}\ell+Y)
+(\Delta_{\mathfrak{g}}\ell+R)u^{2}Y.
\end{aligned}
\right.
\end{equation}
\end{theorem}

{\bf Proof}. Recalling that
\begin{equation*}
|e^{\ell}\Delta_{\mathfrak{g}}v|^{2}=|I_{1}+I_{2}|^{2}=|I_{1}|^{2}+|I_{2}|^{2}+2I_{1}I_{2}.
\end{equation*}
It suffices to compute $2I_{1}I_{2}$. Denote the
terms in the right hand side of $I_{1}$ and
$I_{2}$ by $I_{1}^{d}$ $(d=1,2,3)$ and
$I_{2}^{d}$ $(d=1,2,3)$, respectively. Then
\begin{equation*}
2I_{1}I_{2}=2I_{1}^{1}I_{2}^{1}+2(I_{1}^{1}I_{2}^{2}+I_{1}^{2}I_{2}^{1})+2(I_{1}^{1}I_{2}^{3}
+I_{1}^{3}I_{2}^{1})+2I_{1}^{2}I_{2}^{2}+2(I_{1}^{2}I_{2}^{3}+I_{1}^{3}I_{2}^{2})+2I_{1}^{3}I_{2}^{3}.
\end{equation*}
By \eqref{2.1} and \eqref{2.2}, we compute
\begin{equation}\label{2.6}
\begin{aligned}
2I_{1}^{1}I_{2}^{2}
=&-2\mathrm{div}_{\mathfrak{g}}\,\big[(Q+I)\nabla_{\mathfrak{g}}u\big]\langle 2\nabla_{\mathfrak{g}}\ell+Y,\nabla_{\mathfrak{g}}u\rangle_{\mathfrak{g}} \\
=&-2\mathrm{div}_{\mathfrak{g}}\,\big[\langle
2\nabla_{\mathfrak{g}}\ell+Y,\nabla_{\mathfrak{g}}u\rangle_{\mathfrak{g}}
(Q+I)\nabla_{\mathfrak{g}}u\big]
+2\langle(\mathcal{D}_{\mathfrak{g}}(2\nabla_{\mathfrak{g}}\ell+Y),\nabla_{\mathfrak{g}}u)_{\mathfrak{g}},(Q+I)\nabla_{\mathfrak{g}}u\rangle_{\mathfrak{g}}\\
&+2\langle(2\nabla_{\mathfrak{g}}\ell+Y,\mathcal{D}_{\mathfrak{g}}\nabla_{\mathfrak{g}}u)_{\mathfrak{g}},(Q+I)\nabla_{\mathfrak{g}}u\rangle_{\mathfrak{g}}.
\end{aligned}
\end{equation}
Similarly,
\begin{equation}\label{2.7}
\begin{aligned}
2I_{1}^{2}I_{2}^{1}
=&-2\mathrm{div}_{\mathfrak{g}}\,\big[\langle
Y,\nabla_{\mathfrak{g}}u\rangle_{\mathfrak{g}}Q\nabla_{\mathfrak{g}}u\big]
+2\langle(\mathcal{D}_{\mathfrak{g}}Y,\nabla_{\mathfrak{g}}u)_{\mathfrak{g}},Q\nabla_{\mathfrak{g}}u\rangle_{\mathfrak{g}}\\
&+2\langle(Y,\mathcal{D}_{\mathfrak{g}}\nabla_{\mathfrak{g}}u)_{\mathfrak{g}},Q\nabla_{\mathfrak{g}}u\rangle_{\mathfrak{g}}.
\end{aligned}
\end{equation}
Put
\begin{equation}\label{2.8}
F=\langle(2\nabla_{\mathfrak{g}}\ell+Y,\mathcal{D}_{\mathfrak{g}}\nabla_{\mathfrak{g}}u)_{\mathfrak{g}},(Q+I)\nabla_{\mathfrak{g}}u\rangle_{\mathfrak{g}}
+\langle(Y,\mathcal{D}_{\mathfrak{g}}\nabla_{\mathfrak{g}}u)_{\mathfrak{g}},Q\nabla_{\mathfrak{g}}u\rangle_{\mathfrak{g}}.
\end{equation}
Next, we have
\begin{equation}\label{2.9}
\begin{aligned}
&2(I_{1}^{1}I_{2}^{3}+I_{1}^{3}I_{2}^{1})\\
&=-2\mathrm{div}_{\mathfrak{g}}\,\big[(Q+I)\nabla_{\mathfrak{g}}u\big](\Delta_{\mathfrak{g}}\ell+R)u
-2\mathrm{div}_{\mathfrak{g}}\,(Q\nabla_{\mathfrak{g}}u)(|\nabla_{\mathfrak{g}}\ell|_{\mathfrak{g}}^{2}+R)u\\
&=-2\mathrm{div}_{\mathfrak{g}}\,\big[(\Delta_{\mathfrak{g}}\ell+R)u(Q+I)\nabla_{\mathfrak{g}}u\big]
-2\mathrm{div}_{\mathfrak{g}}\,\big[(|\nabla_{\mathfrak{g}}\ell|_{\mathfrak{g}}^{2}+R)u(Q\nabla_{\mathfrak{g}}u)\big]\\
&\q+2\langle\nabla_{\mathfrak{g}}(\Delta_{\mathfrak{g}}\ell+R),(Q+I)\nabla_{\mathfrak{g}}u\rangle_{\mathfrak{g}}u
+2\langle\nabla_{\mathfrak{g}}(|\nabla_{\mathfrak{g}}\ell|_{\mathfrak{g}}^{2}+R),Q\nabla_{\mathfrak{g}}u\rangle_{\mathfrak{g}}u\\
&\q+2(\Delta_{\mathfrak{g}}\ell+R)\langle\nabla_{\mathfrak{g}}u,(Q+I)\nabla_{\mathfrak{g}}u\rangle_{\mathfrak{g}}
+2(|\nabla_{\mathfrak{g}}\ell|_{\mathfrak{g}}^{2}+R)\langle\nabla_{\mathfrak{g}}u,Q\nabla_{\mathfrak{g}}u\rangle_{\mathfrak{g}},
\end{aligned}
\end{equation}
and
\begin{equation}\label{2.10}
\begin{aligned}
& 2(I_{1}^{2}I_{2}^{3}+I_{1}^{3}I_{2}^{2})\\
& =-2\langle Y,\nabla_{\mathfrak{g}}u\rangle_{\mathfrak{g}}(\Delta_{\mathfrak{g}}\ell+R)u-2\langle 2\nabla_{\mathfrak{g}}\ell+Y,\nabla_{\mathfrak{g}}u\rangle_{\mathfrak{g}}(|\nabla_{\mathfrak{g}}\ell|_{\mathfrak{g}}^{2}+R)u\\
&=-\mathrm{div}_{\mathfrak{g}}\,\big[(\Delta_{\mathfrak{g}}\ell+R)u^{2}Y\big]
-\mathrm{div}_{\mathfrak{g}}\,\big[(|\nabla_{\mathfrak{g}}\ell|_{\mathfrak{g}}^{2}+R)u^{2}(2\nabla_{\mathfrak{g}}\ell+Y)\big]\\
&\q
+\big[(\Delta_{\mathfrak{g}}\ell+R)\mathrm{div}_{\mathfrak{g}}\,Y
+\langle\nabla_{\mathfrak{g}}(\Delta_{\mathfrak{g}}\ell+R),Y\rangle_{\mathfrak{g}}\big]u^{2}\\
&\q
+\big[(|\nabla_{\mathfrak{g}}\ell|_{\mathfrak{g}}^{2}+R)\mathrm{div}_{\mathfrak{g}}\,(2\nabla_{\mathfrak{g}}\ell+Y)
+\langle\nabla_{\mathfrak{g}}(|\nabla_{\mathfrak{g}}\ell|_{\mathfrak{g}}^{2}+R),2\nabla_{\mathfrak{g}}\ell+Y\rangle_{\mathfrak{g}}\big]u^{2}.
\end{aligned}
\end{equation}
We also have
\begin{equation}\label{2.11}
\begin{aligned}
&2I_{1}^{1}I_{2}^{1}+2I_{1}^{2}I_{2}^{2}+2I_{1}^{3}I_{2}^{3}
\\ & = -2\mathrm{div}_{\mathfrak{g}}\,\big[(Q+I)\nabla_{\mathfrak{g}}u\big]\mathrm{div}_{\mathfrak{g}}\,(Q\nabla_{\mathfrak{g}}u)-2\langle Y,\nabla_{\mathfrak{g}}u\rangle_{\mathfrak{g}}\langle 2\nabla_{\mathfrak{g}}\ell+Y,\nabla_{\mathfrak{g}}u\rangle_{\mathfrak{g}}\\
&\q
-2(|\nabla_{\mathfrak{g}}\ell|_{\mathfrak{g}}^{2}+R)(\Delta_{\mathfrak{g}}\ell+R)u^{2}.
\end{aligned}
\end{equation}
\par
Finally, combining \eqref{2.6}--\eqref{2.11}, we
get the desired result immediately.
\endpf

\ms

In the rest of this paper, in the case of
dimension $n\geq3$, we will always use a special
case of Theorem \ref{thm2.1} when $Q=0$ and
$Y=0$  as the following corollary.

\begin{corollary}\label{cor2.1}
We have the following pointwise identity
\begin{equation}\label{2.12}
|e^{\ell}\Delta_{\mathfrak{g}}v|^{2}+\mathrm{div}_{\mathfrak{g}}\,\widetilde{V}=|I_{1}|^{2}+|I_{2}|^{2}+\widetilde{B_{1}}u^{2}+2\widetilde{B_{2}}u+2\widetilde{B_{3}},
\end{equation}
where
\begin{equation}\label{2.13}
\left\{
\begin{aligned}
I_{1}&=\Delta_{\mathfrak{g}}u+(|\nabla_{\mathfrak{g}}\ell|_{\mathfrak{g}}^{2}+R)u,\\
I_{2}&=-2\langle\nabla_{\mathfrak{g}}\ell,\nabla_{\mathfrak{g}}u\rangle_{\mathfrak{g}}-(\Delta_{\mathfrak{g}}\ell+R)u,\\
\end{aligned}
\right.
\end{equation}
and
\begin{equation}\label{2.14}
\left\{
\begin{array}{ll}\ds
\widetilde{B_{1}}=
2\langle\nabla_{\mathfrak{g}}(|\nabla_{\mathfrak{g}}\ell|_{\mathfrak{g}}^{2}+R),\nabla_{\mathfrak{g}}\ell\rangle_{\mathfrak{g}}
-2(|\nabla_{\mathfrak{g}}\ell|_{\mathfrak{g}}^{2}+R)R,\\
\ns\ds \widetilde{B_{2}}= \langle\nabla_{\mathfrak{g}}(\Delta_{\mathfrak{g}}\ell+R),\nabla_{\mathfrak{g}}u\rangle_{\mathfrak{g}},\\
\ns\ds \widetilde{B_{3}}= 2\langle(\mathcal{D}_{\mathfrak{g}}\nabla_{\mathfrak{g}}\ell,\nabla_{\mathfrak{g}}u)_{\mathfrak{g}},\nabla_{\mathfrak{g}}u\rangle_{\mathfrak{g}},\\
\ns\ds \widetilde{V}=
2\big[2\langle\nabla_{\mathfrak{g}}\ell,\nabla_{\mathfrak{g}}u\rangle_{\mathfrak{g}}+(\Delta_{\mathfrak{g}}\ell+R)u\big]\nabla_{\mathfrak{g}}u
+2\big[(|\nabla_{\mathfrak{g}}\ell|_{\mathfrak{g}}^{2}+R)u^{2}-|\nabla_{\mathfrak{g}}u|^{2}_{\mathfrak{g}}\big]\nabla_{\mathfrak{g}}\ell.
\end{array}
\right.
\end{equation}
\end{corollary}

{\bf Proof}. We choose $Q=0$ and $Y=0$ in
Theorem \ref{thm2.1}. From \eqref{2.8} then
\begin{equation*}
F=2\langle(\nabla_{\mathfrak{g}}\ell,\mathcal{D}_{\mathfrak{g}}\nabla_{\mathfrak{g}}u)_{\mathfrak{g}},\nabla_{\mathfrak{g}}u\rangle_{\mathfrak{g}}
=\mathrm{div}_{\mathfrak{g}}\,(|\nabla_{\mathfrak{g}}u|^{2}_{\mathfrak{g}}\nabla_{\mathfrak{g}}\ell)-|\nabla_{\mathfrak{g}}u|^{2}_{\mathfrak{g}}\Delta_{\mathfrak{g}}\ell.
\end{equation*}
Combining this with \eqref{2.4} and \eqref{2.5},
we get the desired result immediately.
\endpf

Next, we specialize Theorem \ref{thm2.1} to the
case of the Euclidean metric $\mathfrak{e}$, but
for operators with variable coefficients. Let
$A(x) = (a^{jk}(x))_{n \times n}$ be a symmetric
invertible matrix, where $a^{jk} \in
C^1(\mathbb{R}^n)$ and $x \in \mathbb{R}^n$. We
define $\mathfrak{g}(x) = A^{-1}(x)$ as a
Riemannian metric on $\mathbb{R}^n$ and consider
$(\mathbb{R}^n, \mathfrak{g})$. Therefore, we
have the following relation:
\begin{equation*}
\langle X, Y \rangle_{\mathfrak{g}} = \langle
A^{-1}(x)X, Y \rangle \quad \text{for} , X, Y
\in \mathbb{R}_x^n, , x \in \mathbb{R}^n.
\end{equation*}

By using the relations between the operators
under the two metrics (the Euclidean metric
denoted as $\mathfrak{e}$ and the Riemannian
metric denoted as $\mathfrak{g}$), we can derive
the following pointwise identity for
second-order partial differential operators with
variable coefficients on $\mathbb{R}^n$.

For $v\in C^{2}(\mathbb{R}^{n};\mathbb{R})$, let
\begin{equation*}
\mathcal{P}v=\sum_{j,k=1}^{n}(a^{jk}v_{x_{j}})_{x_{k}}.
\end{equation*}
Here $a^{jk}\in
C^{1}(\mathbb{R}^{n};\mathbb{R})$ satisfies
$a^{jk}=a^{kj}$ for $j,k=1,2,\cdots,n$. Fix
$\ell\in C^{3}(\mathbb{R}^{n};\mathbb{R})$. Put
\begin{equation*}
\theta=e^{\ell},\quad u=\theta v,
\end{equation*}
We have
\begin{equation*}
\theta\mathcal{P}v=\sum_{j,k=1}^{n}(a^{jk}u_{x_{j}})_{x_{k}}-2\sum_{j,k=1}^{n}a^{jk}\ell_{x_{j}}u_{x_{k}}
+\sum_{j,k=1}^{n}a^{jk}\ell_{x_{j}}\ell_{x_{k}}u-\sum_{j,k=1}^{n}(a^{jk}\ell_{x_{j}})_{x_{k}}u.
\end{equation*}

In order to have more flexibility, we introduce
a symmetric matrix $(p^{jk}(x))_{n\times n}$
where $p^{jk}\in
C^{1}(\mathbb{R}^{n};\mathbb{C})$. Put
\begin{equation}\label{2.15}
\left\{
\begin{aligned}
I_{1}&=\sum_{j,k=1}^{n}\Big[(a^{jk}+p^{jk})u_{x_{j}}\Big]_{x_{k}}+\sum_{k=1}^{n}b^{k}u_{x_{k}}
+\Big(\sum_{j,k=1}^{n}a^{jk}\ell_{x_{j}}\ell_{x_{k}}+R\Big)u,\\
I_{2}&=-\sum_{j,k=1}^{n}(p^{jk}u_{x_{j}})_{x_{k}}-\sum_{k=1}^{n}\Big(\sum_{j=1}^{n}a^{jk}\ell_{x_{j}}+b^{k}\Big)u_{x_{k}}
-\Big[\sum_{j,k=1}^{n}(a^{jk}\ell_{x_{j}})_{x_{k}}+R\Big]u,
\end{aligned}
\right.
\end{equation}
where $R\in C^{1}(\mathbb{R}^{n};\mathbb{C})$
and $b^{k}\in C^{1}(\mathbb{R}^{n};\mathbb{C})$
for $k=1,2,\cdots,n$. In the following, for
$z\in\mathbb{C}$, we denote by $\overline{z}$
the complex conjugate of $z$.

We have the following pointwise identity.
\begin{corollary}\label{cor2.2}
It holds that
\begin{equation}\label{2.16}
\begin{aligned}
&|\theta\mathcal{P}v|^{2}+\mathrm{div}\,V \\
&  =|I_{1}|^{2}+|I_{2}|^{2}+Bu^{2}+2\sum_{j=1}^{n}h^{j}u_{x_{j}}u+\sum_{j,k=1}^{n}c^{jk}u_{x_{j}}u_{x_{k}}\\
& \q-2\mathrm{Re}~\Big\{\sum_{j,k,r,
s=1}^{n}\Big[(a^{jk}+p^{jk})u_{x_{j}}\Big]_{x_{k}}(\overline{p^{rs}}u_{x_{r}})_{x_{s}}\Big\},
\end{aligned}
\end{equation}
where
\begin{equation}\label{2.17}\left\{
\begin{aligned}
V=&[V^{1},\cdots,V^{k},\cdots,V^{n}],\\
V^{k}=&\mathrm{Re}\,\Bigg\{\sum_{r,s=1}^{n}\Big\{\Big[2(a^{rk}+p^{rk})\Big(\sum_{j=1}^{n}2a^{js}\ell_{x_{j}}+\overline{b^{s}}\Big)
-(a^{rs}+p^{rs})\Big(\sum_{j=1}^{n}2a^{jk}\ell_{x_{j}}+\overline{b^{k}}\Big)\\
&\q +\Big(2b^{s}\overline{p^{kr}}-b^{k}\overline{p^{rs}}\Big)\Big]u_{x_{r}}u_{x_{s}}\Big\}\\
&\q+2\sum_{j=1}^{n}\Big[(a^{jk}+p^{jk})\Big(\sum_{r,s=1}^{n}(a^{rs}\ell_{x_{r}})_{x_{s}}+\overline{R}\Big)+
\overline{p^{jk}}\Big(\sum_{r,s=1}^{n}a^{rs}\ell_{x_{r}}\ell_{x_{s}}+R\Big)\Big]u_{x_{j}}u\\
&\q+\Big[b^{k}\Big(\sum_{r,s=1}^{n}(a^{rs}\ell_{x_{r}})_{x_{s}}+\overline{R}\Big)+\Big(\sum_{r,s=1}^{n}a^{rs}\ell_{x_{r}}\ell_{x_{s}}+R\Big)
\Big(\sum_{j=1}^{n}2a^{jk}\ell_{x_{j}}+\overline{b^{k}}\Big)\Big]u^{2}\Bigg\},
\end{aligned}
\right.
\end{equation}
and
\begin{equation}\label{2.18}
\left\{
\begin{aligned}
B=&\mathrm{Re}\,\Bigg\{\sum_{k=1}^{n}\Big[b^{k}\Big(\sum_{r,s=1}^{n}(a^{rs}\ell_{x_{r}})_{x_{s}}+\overline{R}\Big)
+\Big(\sum_{r,s=1}^{n}a^{rs}\ell_{x_{r}}\ell_{x_{s}}+R\Big)
\Big(\sum_{j=1}^{n}2a^{jk}\ell_{x_{j}}+\overline{b^{k}}\Big)\Big]_{x_{k}}\\
&\q-2\Big(\sum_{j,k=1}^{n}a^{jk}\ell_{x_{j}}\ell_{x_{k}}+R\Big)\Big(\sum_{r,s=1}^{n}(a^{rs}\ell_{x_{r}})_{x_{s}}+\overline{R}\Big)\Bigg\},\\
h^{j}=&\mathrm{Re}\,\Bigg\{\sum_{k=1}^{n}\Big\{(a^{jk}+p^{jk})\Big(\sum_{r,s=1}^{n}(a^{rs}\ell_{x_{r}})_{x_{s}}+\overline{R}\Big)_{x_{k}}
+\overline{p^{jk}}\Big(\sum_{r,s=1}^{n}a^{rs}\ell_{x_{r}}\ell_{x_{s}}+R\Big)_{x_{k}}\Big\}\Bigg\},\\
c^{jk}=&\mathrm{Re}\,\Bigg\{\sum_{s=1}^{n}\Big\{2(a^{js}+p^{js})\Big(\sum_{r=1}^{n}2a^{kr}\ell_{x_{r}}+\overline{b^{k}}\Big)_{x_{s}}
-\Big[(a^{jk}+p^{jk})\Big(\sum_{r=1}^{n}2a^{rs}\ell_{x_{r}}+\overline{b^{s}}\Big)\Big]_{x_{s}}\\
&\q +2b_{x_{s}}^{k}\overline{p^{js}}-\Big(b^{s}\overline{p^{jk}}\Big)_{x_{s}}\Big\}-2b^{k}\Big(\sum_{r=1}^{n}2a^{jr}\ell_{x_{r}}+\overline{b^{j}}\Big)\\
&\q
+2\Big[(a^{jk}+p^{jk})\Big(\sum_{r,s=1}^{n}(a^{rs}\ell_{x_{r}})_{x_{s}}+\overline{R}\Big)+
\overline{p^{jk}}\Big(\sum_{r,s=1}^{n}a^{rs}\ell_{x_{r}}\ell_{x_{s}}+R\Big)\Big]\Bigg\}.
\end{aligned}
\right.
\end{equation}
\end{corollary}
\begin{remark}
When the symmetric matrix
$A(x)=(a^{jk}(x))_{n\times n}$ is not
invertible, the result still holds. However, the
weighted identity cannot be obtained directly
from Theorem \ref{thm2.1}. In this situation, we
can use an analogous argument as in the proof of
Theorem \ref{thm2.1} to derive the same result.
\end{remark}
\begin{remark}
Corollary \ref{cor2.2} is a generalization of
the fundamental weighted identity presented in
\cite[Theorem 1.1]{FLZ19}. In this corollary,
when we divide $\theta\mathcal{P}v$ into $I_{1}$
? and $I_{2}$, we introduce a first-order
derivative term in $I_{1}$ and a second-order
derivative term in $I_{2}$.
\end{remark}

\section{Proof of Theorem \ref{thm1.1}}

In this section we prove Theorem \ref{thm1.1}.

{\bf Proof of Theorem \ref{thm1.1}}. It is
suffice to prove the equality \eqref{1.2} for
$v\in C^2(M)$. We divide the proof into four
steps.

\ms

{\bf Step 1}. In this step, we introduce the
weight function $\ell$.

For a positive function $c\in C^2(M)$, we have
\begin{equation*}
c^{\frac{n+2}{4}}(-\Delta_{g}+X+q)u=(-\Delta_{c^{-1}g}+cX+q_{c})(c^{\frac{n-2}{4}}u),
\end{equation*}
where
$q_{c}=c^{\frac{n+2}{4}}\Delta_{g}(c^{-\frac{n-2}{4}})-\frac{n-2}{4}Xc+cq$.
Hence, the equality \eqref{1.2} is invariant
under a conformal change of metrics.
Consequently, we only need to handle the case
that the limiting Carleman weight $\varphi$ is a
distance function, i.e.,
\begin{equation}\label{9.20-eq8}
    |\nabla_{g}\varphi|_{g}^{2}=1.
\end{equation}
Indeed, we can replace $g$ by the conformal
metric
$\overline{g}=|\nabla_{g}\varphi|_{g}^{2}g$ to
have
$|\nabla_{\overline{g}}\varphi|_{\overline{g}}^{2}=1$,
then $\varphi$ is a distance function on
$(M,\overline{g})$.

Next, we choose the weight function
\begin{equation}\label{9.20-eq9}
\ell=\tau
\(\varphi+\frac{\varepsilon}{2}\varphi^{2}\),
\end{equation}
where $\varepsilon$ is a small parameter which
will be fixed later.

\ms

{\bf Step 2}. In this step, we apply Corollary
\ref{cor2.1} with   $R=0$ and $\ell$ given by
\eqref{9.20-eq9}.

For $u=e^{\ell}v$, by choosing $R=0$ in
Corollary \ref{cor2.1}, we obtain that
\begin{equation}\label{3.1}
\begin{array}{ll}\ds
|e^{\ell}\Delta_{g}v|^{2}+\mathrm{div}_{g}\,\widetilde{V}
\\ \ns\ds=
|I_{1}|^{2}+|I_{2}|^{2}+4\mathrm{Hess}_{g}\,\ell(\nabla_{g}\ell,\nabla_{g}\ell)u^{2}
+2\langle\nabla_{g}(\Delta_{g}\ell),\nabla_{g}u\rangle_{g}u
+4\langle(\mathcal{D}_{g}\nabla_{g}\ell,\nabla_{g}u)_{g},\nabla_{g}u\rangle_{g},
\end{array}
\end{equation}
where
\begin{equation}\label{3.2}
\widetilde{V}=2\big[2\langle\nabla_{g}\ell,\nabla_{g}u\rangle_{g}+(\Delta_{g}\ell)u\big]
\nabla_{g}u+2(|\nabla_{g}\ell|_{g}^{2}u^{2}-|\nabla_{g}u|^{2})\nabla_{g}\ell.
\end{equation}

From \eqref{9.20-eq9}, we get that
$$
\nabla_{g}\ell=\tau(1+\varepsilon\varphi)\nabla_{g}\varphi.
$$
By \eqref{9.20-eq8} and \eqref{9.20-eq9}, we see
that
$$
\Delta_{g}\ell=\tau\big[\varepsilon|\nabla_{g}\varphi|_{g}^{2}+(1+\varepsilon\varphi)\Delta_{g}\varphi\big]
=\tau\big[\varepsilon+(1+\varepsilon\varphi)\Delta_{g}\varphi\big].
$$
By \eqref{9.20-eq8} again, we obtain
\begin{equation*}
\mathrm{Hess}_{g}\,\varphi(\nabla_{g}\varphi,\nabla_{g}\varphi)
=\frac{1}{2}\nabla_{g}\varphi\Big(\langle\nabla_{g}\varphi,\nabla_{g}\varphi\rangle_{g}\Big)=0,
\end{equation*}
which, together with \eqref{9.20-eq9}, implies
$$
\begin{aligned}
\mathrm{Hess}_{g}\,\ell(\nabla_{g}\ell,\nabla_{g}\ell)&=\tau^{3}\varepsilon(1+\varepsilon\varphi)^{2}
|\nabla_{g}\varphi|_{g}^{4}+\tau^{3}(1+\varepsilon\varphi)^{3}\mathrm{Hess}_{g}\,\varphi(\nabla_{g}\varphi,
\nabla_{g}\varphi)\\
&=\tau^{3}\varepsilon(1+\varepsilon\varphi)^{2}.
\end{aligned}
$$
By choosing
$\varepsilon<\frac{1}{2\|\varphi\|_{L^\infty(M)}}$,
we have
\begin{equation}\label{3.3-1}
\frac{1}{2}\leq
1+\varepsilon\varphi\leq\frac{3}{2}.
\end{equation}
By Cauchy-Schwarz inequality, we obtain that
\begin{equation}\label{3.3}
|2\langle\nabla_{g}(\Delta_{g}\ell),\nabla_{g}u\rangle_{g}u|\leq
C(\tau^{2}u^{2}+|\nabla_{g}u|_{g}^{2})
\end{equation}
and that
\begin{equation}\label{3.4}
|4\langle(\mathcal{D}_{g}\nabla_{g}\ell,\nabla_{g}u)_{g},\nabla_{g}u\rangle_{g}|
\leq C\tau|\nabla_{g}u|_{g}^{2}.
\end{equation}
Combining \eqref{3.1}, \eqref{3.3}  and
\eqref{3.4}, we conclude that
\begin{equation}\label{3.5}
\big|e^{\tau\varepsilon\varphi^{2}/2}e^{\tau\varphi}\Delta_{g}v
\big|^{2}+\mathrm{div}_{g}\,\widetilde{V} \geq
|I_{1}|^{2}+|I_{2}|^{2}+C_{1}\tau^{3}\varepsilon
u^{2}
-C_{2}\tau^{2}u^{2}-C_{3}|\nabla_{g}u|_{g}^{2}.
\end{equation}

For $I_{2}=-2\langle
\nabla_{g}\ell,\nabla_{g}u\rangle_{g}-(\Delta_{g}\ell)
u$, by using the inequality
$(2\alpha+\beta)^{2}\geq2\alpha^{2}-\beta^{2}$,
we find that
\begin{equation}\label{9.20-eq10}
|I_{2}|^{2}\geq
C_{4}\tau^{2}|\nabla_{g}u|_{g}^{2}-C_{5}\tau^{2}u^{2}.
\end{equation}
For $\tau\geq 1$,  we get from \eqref{9.20-eq10}
that
\begin{equation}\label{3.6}
|I_{2}|^{2}\geq
\tau^{-1}\varepsilon|I_{2}|^{2}\geq
\tau\varepsilon
\big(C_{4}|\nabla_{g}u|_{g}^{2}-C_{5}u^{2}\big).
\end{equation}
From \eqref{3.5} and \eqref{3.6}, we conclude
that, for $\tau\geq 1$,
\begin{equation}\label{3.7}
\big|e^{\tau\varepsilon\varphi^{2}/2}e^{\tau\varphi}\Delta_{g}v\big|^{2}+\mathrm{div}_{g}\,\widetilde{V}
\geq
\big(C_{1}\tau^{3}\varepsilon-C_{2}\tau^{2}-C_{5}\tau\varepsilon\big)
u^{2}
+\big(C_{4}\tau\varepsilon-C_{3}\big)|\nabla_{g}u|_{g}^{2}.
\end{equation}

Next, a direct computation yields
\begin{equation}\label{3.7-1}
\begin{aligned}
e^{\ell}(-\Delta_{g}+X+q)v & =-e^{\ell}\Delta_{g}v+e^{\ell}(X+q)v \\
&=-e^{\ell}\Delta_{g}v+(X+q)u-\tau(1+\varepsilon\varphi)X(\varphi
u).
\end{aligned}
\end{equation}
Recalling $u=e^{\ell}v$, we get that
\begin{equation}\label{3.8}
|e^{\ell}(X+q)v|^{2}\leq
C\big(\tau^{2}u^{2}+|\nabla_{g}u|_{g}^{2}\big).
\end{equation}

Using \eqref{3.7}--\eqref{3.8} and the
inequality
\begin{equation}\label{3.8-1}
|e^{\ell}(-\Delta_{g}+X+q)v|^{2} \geq
\frac{1}{2}|e^{\ell}\Delta_{g}v|^{2}
-|e^{\ell}(X+q)v|^{2},
\end{equation}
we obtain that
\begin{equation}\label{3.9}
\begin{array}{ll}\ds
\int_{M}|e^{\tau\varepsilon\varphi^{2}/2}e^{\tau\varphi}(-\Delta_{g}+X+q)v|^{2}\,dV_{g} + \int_{M}\mathrm{div}_{g}\,\widetilde{V}\,dV_{g}\\
\ns\ds\geq
\int_{M}\big[(C_{1}\tau^{3}\varepsilon-C_{2}\tau^{2}-C_{5}\tau\varepsilon)
u^{2}
+(C_{4}\tau\varepsilon-C_{3})|\nabla_{g}u|_{g}^{2}
\big]\,dV_{g}.
\end{array}
\end{equation}

\ms

{\bf Step 3}. In this step, we deal with the
term $\int_{M}\mathrm{div}_{g}\,\widetilde{V}
dV_{g}$ in the left hand side of \eqref{3.9}.

On the boundary $\pa M$,
\begin{equation*}
\nabla_{g}u=\nabla_{\perp}u+\nabla_{\parallel}u,
\end{equation*}
where
$\nabla_{\perp}u=\langle\nabla_{g}u,\nu\rangle_{g}\nu$
and
$\nabla_{\parallel}u=\nabla_{g}u-\nabla_{\perp}u$.
Recall the expression of $\widetilde{V}$ in
\eqref{3.2}, we have
\begin{equation}\label{3.10}
\begin{aligned}
&\int_{M}\mathrm{div}_{g}\,\widetilde{V}\,dV_{g}\\&=  \int_{\pa M}\langle\mathrm{div}_{g}\,\widetilde{V},\nu\rangle_{g}\,dS_g\\
&=2\int_{\pa M}\Big[(2\langle\nabla_{g}\ell,\nabla_{g}u\rangle_{g}+u\Delta_{g}\ell) \langle\nabla_{g}u,\nu\rangle_{g}+(|\nabla_{g}\ell|_{g}^{2}u^{2}-|\nabla_{g}u|^{2})\langle\nabla_{g}\ell,\nu\rangle_{g}\Big]\,  dS_g\\
&=\int_{\pa M}
\[2\tau(1+\varepsilon\varphi)\langle\nabla_{\perp}\varphi,\nu\rangle_{g}|\nabla_{\perp}u|_{g}^{2}
+4\tau(1+\varepsilon\varphi)\langle\nabla_{\parallel}\varphi,\nabla_{\parallel}u\rangle_{g}\langle\nabla_{\perp}u,\nu\rangle_{g}\\
&\qq\qq
+2\tau\big[\varepsilon+(1+\varepsilon\varphi)\Delta_{g}\varphi\big](\langle\nabla_{\perp}u,\nu\rangle_{g}
+\langle\nabla_{\parallel}u,\nu\rangle_{g})u\\
&\qq\qq
+2\tau^{3}(1+\varepsilon\varphi)^{3}\langle\nabla_{\perp}\varphi,\nu\rangle_{g}u^{2}
-2\tau(1+\varepsilon\varphi)\langle\nabla_{\perp}\varphi,\nu\rangle_{g}|\nabla_{\parallel}u|_{g}^{2} \] \,dS_g\\
&\leq C\int_{\pa
M}\big(\tau|\langle\nabla_{\perp}\varphi,\nu\rangle_{g}|\,|\nabla_{\perp}u|_{g}^{2}+
\tau|\nabla_{\parallel}u|_{g}|\nabla_{\perp}u|_{g}+\tau|\nabla_{\perp}u|_{g}|u|+\tau|\nabla_{\parallel}u|_{g}^{2}
+\tau^{3}|u|^{2}\big)\, dS_g.
\end{aligned}
\end{equation}
From the equality
\begin{equation}\label{9.20-eq11}
\nabla_{g}u=e^{\ell}\big[\nabla_{g}v+\tau(1+\varepsilon\varphi)v\nabla_{g}\varphi\big]
=e^{\ell}\nabla_{g}v+\tau(1+\varepsilon\varphi)u\nabla_{g}\varphi,
\end{equation}
we obtain that
\begin{eqnarray}
|\nabla_{\parallel}u|_{g}\3n &\leq
Ce^{\ell}\big(|\nabla_{\parallel}v|_{g}+\tau|v|\big),\label{3.11}
\\ [2mm] |\nabla_{\perp}u|_{g}\3n &\leq
Ce^{\ell}\big(|\nabla_{\perp}v|_{g}+\tau|v|\big).\label{3.12}
\end{eqnarray}
By \eqref{3.11}, \eqref{3.12} and the
Cauchy-Schwarz inequality, we get that
\begin{equation}\label{3.14}
\begin{aligned}
\tau|\nabla_{\parallel}u|_{g}|\nabla_{\perp}u|_{g}&\leq
C\tau e^{2\ell}(|\nabla_{\parallel}v|_{g}+\tau|v|)(|\nabla_{\perp}v|_{g}+\tau|v|)\\
&\leq
Ce^{2\ell}(\tau|\nabla_{\parallel}v|_{g}|\nabla_{\perp}v|_{g}+\tau^{2}|\nabla_{\parallel}v|_{g}|v|
+\tau^{2}|\nabla_{\perp}v|_{g}|v|+\tau^{3}|v|^{2})\\
&\leq
Ce^{2\ell}(\tau|\nabla_{\parallel}v|_{g}|\nabla_{\perp}v|_{g}+\tau^{2}|\nabla_{\perp}v|_{g}|v|
+\tau|\nabla_{\parallel}v|_{g}^{2}+\tau^{3}|v|^{2}).
\end{aligned}
\end{equation}
It follows from \eqref{3.10}--\eqref{3.14} that
\begin{equation}\label{3.15}
\begin{aligned}
&\int_{M}\mathrm{div}_{g}\,\widetilde{V}\,dV_{g} \\
& \leq   C\int_{\pa
M}e^{\tau\varepsilon\varphi^{2}}\(\tau
e^{2\tau\varphi}|\langle\nabla_{\perp}\varphi,\nu\rangle_{g}|\,|\nabla_{\perp}v|_{g}^{2}+
\tau e^{2\tau\varphi}|\nabla_{\parallel}v|_{g}|\nabla_{\perp}v|_{g}\\
&\qq \qq \qq\;\,
+\tau^{2}e^{2\tau\varphi}|\nabla_{\perp}v|_{g}|v|+\tau
e^{2\tau\varphi}|\nabla_{\parallel}v|_{g}^{2}
+\tau^{3}|e^{\tau\varphi}v|^{2}\)\,dS_g.
\end{aligned}
\end{equation}

\ms

{\bf Step 4}. In this step, we complete the
proof.

Combining \eqref{3.9} and \eqref{3.15}, we  get
that
\begin{equation}\label{3.16}
\begin{aligned}
&
\int_{M}e^{\tau\varepsilon\varphi^{2}}\big|e^{\tau\varphi}(-\Delta_{g}+X+q)v\big|^{2}\,dV_{g}
\\ &  \q +\int_{\partial M}e^{\tau\varepsilon\varphi^{2}}\Big[\tau e^{2\tau\varphi}|\langle\nabla_{\perp}\varphi,\nu\rangle_{g}|\,|\nabla_{\perp}v|_{g}^{2}+
\tau e^{2\tau\varphi}|\nabla_{\parallel}v|_{g}|\nabla_{\perp}v|_{g}\\
& \qq\qq\qq\;
+\tau^{2}e^{2\tau\varphi}|\nabla_{\perp}v|_{g}|v|+\tau|e^{\tau\varphi}\nabla_{\parallel}v|_{g}^{2}
+\tau^{3}|e^{\tau\varphi}v|^{2}\Big]\,dS_{g}\\
&\geq
\int_{M}\big(C_{1}\tau^{3}\varepsilon-C_{2}\tau^{2}-C_{5}\tau\varepsilon
\big) u^{2}
+\big(C_{4}\tau\varepsilon-C_{3}\big)|\nabla_{g}u|_{g}^{2}\,dV_{g}.
\end{aligned}
\end{equation}

By \eqref{9.20-eq11} again, we find that
\begin{equation}\label{3.13}
    \begin{aligned}
        \tau^{2}|e^{\ell}v|^{2}+
        |e^{\ell}\nabla_{g}v|_{g}^{2}&=\tau^{2}|u|^{2}+|\nabla_{g}u
        -\tau(1+\varepsilon\varphi)u\nabla_{g}\varphi|_{g}^{2}\\
        &\leq C(\tau^{2}|u|^{2}+|\nabla_{g}u|_{g}^{2}).
    \end{aligned}
\end{equation}

Let $$C_6\deq\max\Big\{\frac{C_2}{C_1}+1,
\frac{C_3}{C_4}+1\Big\}$$ and $$\tau_0 \deq
\max\Big\{2\(\frac{C_5C_6}{C_1C_6-C_2}\),1\Big\}.$$
By taking $\varepsilon=C_6 \tau^{-1}$, for any
$\tau>\tau_1$, we obtain from \eqref{3.13} that
\begin{equation}\label{3.17}
\begin{aligned}
&\big(C_{1}\tau^{3}\varepsilon-C_{2}\tau^{2}-C_{5}\tau\varepsilon\big)
u^{2}
+\big(C_{4}\tau\varepsilon-C_{3}\big)|\nabla_{g}u|_{g}^{2}\\
&\geq  \big(C_{1}C_6\tau^{2}
-C_{2}\tau^{2}-C_{5}C_6\big) u^{2}
+\big(C_{4}C_6-C_{3}\big)|\nabla_{g}u|_{g}^{2}
\\
&\geq C\big(\tau^{2}|u|^{2}+|\nabla_{g}u|_{g}^{2}\big)\\
&\geq
Ce^{\tau\varepsilon\varphi^{2}}\big(\tau^{2}|e^{\tau\varphi}v|^{2}+|e^{\tau\varphi}
\nabla_{g}v|_{g}^{2}\big).
\end{aligned}
\end{equation}
By \eqref{3.16} and \eqref{3.17},  we derive
that
\begin{equation}\label{3.18}
\begin{aligned}
&
\int_{M}e^{\tau\varepsilon\varphi^{2}}\big|e^{\tau\varphi}(-\Delta_{g}+X+q)v\big|^{2}\,dV_{g}
\\ &  \q +\int_{\partial M}e^{\tau\varepsilon\varphi^{2}}\Big(\tau e^{2\tau\varphi}|\langle\nabla_{\perp}\varphi,\nu\rangle_{g}|\,|\nabla_{\perp}v|_{g}^{2}+
\tau e^{2\tau\varphi}|\nabla_{\parallel}v|_{g}|\nabla_{\perp}v|_{g}\\
& \qq\qq\qq\;
+\tau^{2}e^{2\tau\varphi}|\nabla_{\perp}v|_{g}|v|+\tau|e^{\tau\varphi}\nabla_{\parallel}v|_{g}^{2}
+\tau^{3}|e^{\tau\varphi}v|^{2}\Big)\,dS_{g}\\
&\geq C
\int_{M}e^{\tau\varepsilon\varphi^{2}}\big(\tau^{2}|e^{\tau\varphi}v|^{2}+|e^{\tau\varphi}
\nabla_{g}v|_{g}^{2}\big)\,dV_{g}.
\end{aligned}
\end{equation}
By the choise of $\e$, we know  that
$$0\leq \tau\varepsilon\varphi^{2}\leq C_{6} \|\varphi\|_{L^\infty(M)}^2. $$
This, together with \eqref{3.18}, implies the
equality \eqref{1.2} immediately.
\endpf

\section{Proof of Theorem \ref{thm1.2}}

In this section we establish the Carleman
estimate \eqref{1.3} on Riemann surfaces. To
begin with, we recall the following result.

\begin{lemma}\cite[Proposition 18.9]{LRLR22}
Let $\Gamma_{*}$ be a nonempty open subset of
$\partial \widetilde{M}$. There exists a
constant $C>0$ such that for any $W\in
H^{1}(\widetilde{M})$,
\begin{equation}\label{4.12}
\|W\|_{L^{2}(\widetilde{M})}^{2}\leq
C\big(\|\nabla
W\|_{L^{2}(\widetilde{M})}^{2}+\|W\|_{L^{2}(\Gamma_{*})}^{2}\big).
\end{equation}
\end{lemma}

We observe that the inequality \eqref{1.3} is
analogous to \eqref{1.2}, and the proof of
Theorem \ref{thm1.1} is independent of the
dimension. However, since the weight functions
have some critical points in $\widetilde{M}$, we
need a completely different proof.

\ms

{\bf Proof of Theorem \ref{thm1.2}}.  We divide
the proof into six steps.

\ms

{\bf Step 1}. In this step, we do some
reductions. As we assume that $\widetilde{M}$ is
simply connected, by Proposition 2.4 in
\cite{MT02}, we can choose $\widetilde{M}$ to be
the closed unit disk
$\{z\in\mathbb{C}:|z|\leq1\}$ and such that the
metric $\widetilde{g}$ is conformal to the
Euclidean metric $\mathfrak{e}$, i.e., there
exists a smooth positive function $\lambda(x)$
such that
$\widetilde{g}=e^{2\lambda}\mathfrak{e}$. Thus
the norm induced by $\widetilde{g}$ is conformal
to the Euclidean norm, and the Laplace-Beltrami
operator with respect to $\widetilde{g}$ is
given by
$\Delta_{\widetilde{g}}=e^{-2\lambda}\Delta$.
Then it is suffice  to prove the estimate
\eqref{1.7} under the Euclidean metric.

Furthermore,   it suffices to establish
\eqref{1.7} for $v\in C^{\infty}(\widetilde{M})$
with $v\big|_{\partial \widetilde{M}}=0$.

\ms

{\bf Step 2}. In this step, we apply Corollary
\ref{cor2.2} with suitable choosing $\ell$,
$(p^{jk})_{1\leq j,k\leq 2}$, $(b^{1},b^{2})$
and $R$.

Let $\ell=\tau\varphi$, where
$\varphi:\widetilde{M}\rightarrow\mathbb{R}$ is
the harmonic Morse function as mentioned in
Theorem \ref{thm1.2}. For $u=e^{\ell}v$, let
\begin{eqnarray}\label{9.20-eq1}
\begin{cases}\ds
(p^{jk})_{1\leq j,k\leq 2}=
\begin{pmatrix}
-1 & \frac{i}{2} \\
\frac{i}{2} & 0
\end{pmatrix} , \\
\ns\ds
(b^{1},b^{2})=(-\ell_{x_{1}},-\ell_{x_{2}})=(-\tau\varphi_{x_{1}},-\tau\varphi_{x_{2}}),\\
\ns\ds R=
-(\tau^{2}\varphi_{x_{2}}^{2}+\tau\varphi_{x_{2}x_{2}})-i(\tau^{2}\varphi_{x_{1}}\varphi_{x_{2}}
+\tau\varphi_{x_{1}x_{2}})
\end{cases}
\end{eqnarray}
in Corollary \ref{cor2.2}. Then we have
\begin{equation}\label{4.0}
\begin{aligned}
&|e^{\tau\varphi}\Delta v|^{2}+\mathrm{div}\,V \\
&  =|I_{1}|^{2}+|I_{2}|^{2}+Bu^{2}+2\sum_{j=1}^{n}h^{j}u_{x_{j}}u+\sum_{j,k=1}^{n}c^{jk}u_{x_{j}}u_{x_{k}}\\
& \q-2\mathrm{Re}~\Big\{\sum_{j,k,r,
s=1}^{n}\Big[(a^{jk}+p^{jk})u_{x_{j}}\Big]_{x_{k}}(\overline{p^{rs}}u_{x_{r}})_{x_{s}}\Big\},
\end{aligned}
\end{equation}
It follows from \eqref{9.20-eq1} that
\begin{equation}\label{4.1}
\begin{aligned}
&-2\Re\Big\{\sum_{j,k,r,s=1}^{2}\Big[(a^{jk}+p^{jk})u_{x_{j}}\Big]_{x_{k}}\Big(\overline{p^{rs}}u_{x_{r}}\Big)_{x_{s}}\Big\}\\
&=-2(u_{x_{1}x_{2}}u_{x_{1}x_{2}}-u_{x_{1}x_{1}}u_{x_{2}x_{2}})\\
&=2\Big[(u_{x_{1}}u_{x_{2}x_{2}})_{x_{1}}-(u_{x_{1}}u_{x_{1}x_{2}})_{x_{2}}\Big],
\end{aligned}
\end{equation}
and that
\begin{eqnarray}\label{4.2}
&&B =-|\nabla\ell|^{2}\Delta\ell-(\Delta\ell)^{2}+\nabla\ell\cdot\nabla(|\nabla\ell|^{2}-\Delta\ell)-2(|\nabla\ell|^{2}+\Delta\ell)\mathrm{Re}\,R-2R\overline{R}\nonumber\\
&&\q=2\tau^{3}\sum_{j,k=1}^{2}\varphi_{x_{j}x_{k}}\varphi_{x_{j}}\varphi_{x_{k}}
+2\tau^{2}|\nabla\varphi|^{2}(\tau^{2}\varphi_{x_{2}}^{2}+\tau\varphi_{x_{2}x_{2}})\\
&&\q\q-2\Big[(\tau^{2}\varphi_{x_{2}}^{2}+\tau\varphi_{x_{2}})^{2}+(\tau^{2}\varphi_{x_{1}}\varphi_{x_{2}}+\tau\varphi_{x_{1}x_{2}})^{2}\Big]\nonumber\\
&&\q=-\tau^{2}\sum_{j,k=1}^{2}\varphi_{x_{j}x_{k}}^{2}.\nonumber
\end{eqnarray}
Here and in what follows, we use the fact
$\Delta\varphi=\varphi_{x_{1}x_{1}}+\varphi_{x_{2}x_{2}}=0$
to eliminate some terms.

Since
\begin{equation*}
\begin{cases}
h^{1}=(\mathrm{Im}\,R)_{x_{2}}-(|\nabla\ell|^{2}+\mathrm{Re}\,R)_{x_{1}},\\
h^{2}=(\mathrm{Im}\,R)_{x_{1}}+(\Delta\ell+\mathrm{Re}\,R)_{x_{2}},
\end{cases}
\end{equation*}
we have
\begin{equation}\label{4.3}
\begin{aligned}
2\sum_{j=1}^{2}h^{j}u_{x_{j}}u
=&2u\Big[(-\tau^{2}\varphi_{x_{1}}\varphi_{x_{2}}-\tau\varphi_{x_{1}x_{2}})_{x_{2}}u_{x_{1}}
-(\tau^{2}\varphi_{x_{1}}^{2}-\tau\varphi_{x_{2}x_{2}})_{x_{1}}u_{x_{1}}\\
&\q+(-\tau^{2}\varphi_{x_{1}}\varphi_{x_{2}}-\tau\varphi_{x_{1}x_{2}})_{x_{1}}u_{x_{2}}
-(\tau^{2}\varphi_{x_{2}}^{2}+\tau\varphi_{x_{2}x_{2}})_{x_{2}}u_{x_{2}}\Big]\\
=&-2\tau^{2}\sum_{j,k=1}^{2}\varphi_{x_{j}x_{k}}\varphi_{x_{j}}u_{x_{k}}u.
\end{aligned}
\end{equation}
From \eqref{2.18} and \eqref{9.20-eq1}, we find
that
\begin{equation}\label{9.20-eq2}
c^{jk}=2\big[\ell_{x_{j}x_{k}}+\ell_{x_{j}}\ell_{x_{k}}+a^{jk}\mathrm{Re}\,R+|\nabla\ell|^{2}
\mathrm{Re}\,p^{jk}+2\mathrm{Re}\,(p^{jk}\overline{R})\big].
\end{equation}
Recalling the choice of $p^{jk}$ in
\eqref{9.20-eq1}, we get that
\begin{equation*}
\begin{cases}\ds
c^{11}
\3n&=2\big[\ell_{x_{1}x_{1}}+\ell_{x_{1}}\ell_{x_{1}}+\mathrm{Re}\,R-|\nabla\ell|^{2}
-2\mathrm{Re}\,\overline{R})\big]\\
\ns&\ds
=2\big[\tau\varphi_{x_{1}x_{1}}+\tau^{2}\varphi_{x_{1}}^{2}
+(\tau^{2}\varphi_{x_{2}}^{2}+\tau\varphi_{x_{2}x_{2}})-\tau^{2}(\varphi_{x_{1}}^{2}+\varphi_{x_{2}}^{2})\big]=0,\\
\ns\ds c^{12} \3n&=c^{21}=2\big[\ell_{x_{1}x_{2}}+\ell_{x_{1}}\ell_{x_{2}}+2\mathrm{Re}\,(\frac{i}{2}\overline{R})\big]\\
\ns&\ds
=2\big[\tau\varphi_{x_{1}x_{2}}+\tau^{2}\varphi_{x_{1}}\varphi_{x_{2}}
-(\tau^{2}\varphi_{x_{1}x_{2}}+\tau\varphi_{x_{1}x_{2}})\big]=0,\\
\ns\ds c^{22} \3n&=2\big[\ell_{x_{2}x_{2}}+\ell_{x_{2}}\ell_{x_{2}}+\mathrm{Re}\,R\big]\\
\ns&\ds
=2\big[\tau\varphi_{x_{2}x_{2}}+\tau^{2}\varphi_{x_{2}}^{2}-(\tau^{2}\varphi_{x_{2}}^{2}
+\tau\varphi_{x_{2}x_{2}})\big]=0.
\end{cases}
\end{equation*}
Consequently,
\begin{equation}\label{4.4}
\sum_{j,k=1}^{2}c^{jk}u_{x_{j}}u_{x_{k}}=0.
\end{equation}

Next,
\begin{equation}\label{4.5}
\begin{aligned}
&-\tau^{2}\sum_{j,k=1}^{2}\varphi_{x_{j}x_{k}}^{2}|u|^{2}\\
=&-\tau^{2}\sum_{j,k=1}^{2}(\varphi_{x_{j}x_{k}}\varphi_{x_{j}}|u|^{2})_{x_{k}}
+\tau^{2}\sum_{j,k=1}^{2}(\varphi_{x_{k}x_{k}})_{x_{j}}\varphi_{x_{j}}|u|^{2}
+2\tau^{2}\sum_{j,k=1}^{2}\varphi_{x_{j}x_{k}}\varphi_{x_{j}}v_{x_{k}}u\\
=&-\tau^{2}\sum_{j,k=1}^{2}(\varphi_{x_{j}x_{k}}\varphi_{x_{j}}|u|^{2})_{x_{k}}
+2\tau^{2}\sum_{j,k=1}^{2}\varphi_{x_{j}x_{k}}\varphi_{x_{j}}u_{x_{k}}u.
\end{aligned}
\end{equation}

Integrating \eqref{4.0} on $\widetilde{M}$, from
\eqref{4.1}--\eqref{4.5}, we obtain that
\begin{equation}\label{4.7}
\begin{aligned}
& \int_{\wt M}|e^{\tau\varphi}\Delta v|^{2}\,dx + \int_{\wt M}\mathrm{div}\,V\,dx \\
&=\int_{\wt
M}\big(|I_{1}|^{2}+|I_{2}|^{2}\big)\,dx
+2\int_{\wt
M}\big[(u_{x_{1}}u_{x_{2}x_{2}})_{x_{1}}-(u_{x_{1}}u_{x_{1}x_{2}})_{x_{2}}\big]\,dx
\\
& \q -\tau^{2} \int_{\wt
M}\sum_{j,k=1}^{2}\big(\varphi_{x_{j}x_{k}}\varphi_{x_{j}}|u|^{2}\big)_{x_{k}}\,dx,
\end{aligned}
\end{equation}
where
\begin{equation}\label{4.7-1}\left\{
\begin{aligned}
I_{1}&=u_{x_{2}x_{2}}+iu_{x_{1}x_{2}}-\tau\nabla\varphi\cdot\nabla
u
+\Big[\tau^{2}(\varphi_{x_{1}}^{2}-i\varphi_{x_{1}}\varphi_{x_{2}})
+\tau(-\varphi_{x_{2}x_{2}}-i\varphi_{x_{1}x_{2}})\Big]u,\\
I_{2}&=u_{x_{1}x_{1}}-iu_{x_{1}x_{2}}-\tau\nabla\varphi\cdot\nabla
u
+\Big[\tau^{2}(\varphi_{x_{2}}^{2}+i\varphi_{x_{1}}\varphi_{x_{2}})
+\tau(-\varphi_{x_{1}x_{1}}+i\varphi_{x_{1}x_{2}})\Big]u.
\end{aligned}\right.
\end{equation}

\ms

{\bf Step 3}. In this step, we deal with the
terms concerning the integral on $\partial
\widetilde{M}$ in \eqref{4.7}.

Since $u\big|_{\partial \widetilde{M}}=0$, we
see that $u_{x_{1}} = \partial_{\nu}u\, \nu_{1}
$ on $\partial \widetilde{M}$. Hence, we find
that
\begin{equation}\label{9.20-eq6}
u_{x_{1}x_{2}}=(\partial_{\nu}u)_{x_{2}}\nu_{1}+\partial_{\nu}u(\nu_{1})_{x_{2}}\q
\mbox{ on } \partial \widetilde{M}.
\end{equation}
Similarly, we can obtain
\begin{equation}\label{9.20-eq7}
u_{x_{2}x_{2}}=(\partial_{\nu}u)_{x_{2}}\nu_{2}+\partial_{\nu}u(\nu_{2})_{x_{2}}\q
\mbox{ on } \partial \widetilde{M}.
\end{equation}
It follows from \eqref{9.20-eq6} and
\eqref{9.20-eq7} that
\begin{equation}\label{4.8}
\begin{array}{ll}\ds
2\int_{\wt M}\big[(u_{x_{1}}u_{x_{2}x_{2}})_{x_{1}}-(u_{x_{1}}u_{x_{1}x_{2}})_{x_{2}}\big]\,dx\\
\ns\ds =
2 \int_{\pa\wt M} (u_{x_{1}}u_{x_{2}x_{2}}\nu_{1}-u_{x_{1}}u_{x_{1}x_{2}}\nu_{2})\,dS\\
\ns\ds =2\int_{\pa\wt M}
|\partial_{\nu}u|^{2}\Big[(\nu_{2})_{x_{2}}\nu_{1}^{2}-(\nu_{1})_{x_{2}}\nu_{1}\nu_{2}\Big]\,dS.
\end{array}
\end{equation}
As mentioned in the beginning of the proof, we
can choose the boundary to be the circle $|z|=1$
in $\mathbb{C}$. Then from \eqref{4.8}, we get
that
\begin{equation}\label{4.8-1}
2\int_{\wt
M}\big[(u_{x_{1}}u_{x_{2}x_{2}})_{x_{1}}-(u_{x_{1}}u_{x_{1}x_{2}})_{x_{2}}\big]\,dx
=2\int_{\pa\wt
M}|\partial_{\nu}u|^{2}x_{1}^{2}\,dS.
\end{equation}

Using the boundary condition $u\big|_{\partial
\widetilde{M}}=0$ again, we get that
\begin{equation}\label{4.6}
\begin{aligned}
\int_{\wt M} \mathrm{div}\,V\,dx& = \int_{\pa\wt
M}  V\cdot\nu\,dS
\\ &=\mathrm{Re} \int_{\pa\wt M} \sum_{k,r,s=1}^{2}\Big[2\Big(a^{kr}+p^{kr}-\overline{p^{kr}}\Big)\ell_{x_{s}}
-\Big(a^{rs}+p^{rs}-\overline{p^{rs}}\Big)\ell_{x_{k}}\Big]u_{x_{r}}u_{x_{s}}\nu^{k}\,dS\\
&= \int_{\pa\wt M} \big[2\tau(\nabla\varphi\cdot\nabla u)\partial_{\nu}u-\tau(\nabla\varphi\cdot\nu)|\nabla u|^{2} \big] \,dS\\
&=\tau \int_{\pa\wt M}
\partial_{\nu}\varphi|\partial_{\nu}u|^{2} \,dS.
\end{aligned}
\end{equation}

Combining \eqref{4.7}, \eqref{4.8-1} and
\eqref{4.6}, we find that
\begin{equation}\label{4.9}
\big\|e^{\tau\varphi}\Delta
v\big\|_{L^{2}(\widetilde{M})}^{2}+\tau\int_{\partial
\widetilde{M}}\partial_{\nu}\varphi|\partial_{\nu}u|^{2}\,dS=\big\|I_{1}\big\|_{L^{2}(\widetilde{M})}^{2}+\big\|I_{2}\big\|_{L^{2}(\widetilde{M})}^{2}
+2\int_{\partial
\widetilde{M}}x_{1}^{2}|\partial_{\nu}u|^{2}\,dS.
\end{equation}

Next, we set
\begin{equation*}
\begin{cases}\ds
(p^{jk})_{2\times 2}=
\begin{pmatrix}
-1 & -\frac{i}{2} \\
-\frac{i}{2} & 0
\end{pmatrix} , \\
\ns\ds
(b^{1},b^{2})=(-\ell_{x_{1}},-\ell_{x_{2}})=(-\tau\varphi_{x_{1}},-\tau\varphi_{x_{2}}),\\
\ns\ds    R=
-(\tau^{2}\varphi_{x_{2}}^{2}+\tau\varphi_{x_{2}x_{2}})+i(\tau^{2}\varphi_{x_{1}}\varphi_{x_{2}}
+\tau\varphi_{x_{1}x_{2}})
\end{cases}
\end{equation*}
in Corollary \ref{cor2.2}. Similar to the proof
of \eqref{4.9}, we can obtain that
\begin{equation}\label{4.10}
\big\|e^{\tau\varphi}\Delta
v\big\|_{L^{2}(\widetilde{M})}^{2}+\tau\int_{\partial
\widetilde{M}}\partial_{\nu}\varphi|\partial_{\nu}u|^{2}\,dS=\big\|I_{3}\big\|_{L^{2}(\widetilde{M})}^{2}+\big\|I_{4}\big\|_{L^{2}(\widetilde{M})}^{2}
+2\int_{\partial
\widetilde{M}}x_{2}^{2}|\partial_{\nu}u|^{2}\,dS,
\end{equation}
where
\begin{equation*}\left\{
\begin{aligned}
I_{3}&=u_{x_{2}x_{2}}-iu_{x_{1}x_{2}}-\tau\nabla\varphi\cdot\nabla
u
+\Big[\tau^{2}(\varphi_{x_{1}}^{2}+i\varphi_{x_{1}}\varphi_{x_{2}})
+\tau(-\varphi_{x_{2}x_{2}}+i\varphi_{x_{1}x_{2}})\Big]u,\\
I_{4}&=u_{x_{1}x_{1}}+iu_{x_{1}x_{2}}-\tau\nabla\varphi\cdot\nabla
u
+\Big[\tau^{2}(\varphi_{x_{2}}^{2}-i\varphi_{x_{1}}\varphi_{x_{2}})
+\tau(-\varphi_{x_{1}x_{1}}-i\varphi_{x_{1}x_{2}})\Big]u.
\end{aligned}\right.
\end{equation*}
From \eqref{4.7-1}, since
$I_{3}=\overline{I_{1}}$ and
$I_{4}=\overline{I_{2}}$, by adding \eqref{4.10}
to \eqref{4.9}, we conclude that
\begin{equation}\label{4.11}
\big\|e^{\tau\varphi}\Delta
v\big\|_{L^{2}(\widetilde{M})}^{2}+\tau\int_{\partial
\widetilde{M}}\partial_{\nu}\varphi|\partial_{\nu}u|^{2}\,dS
=\big\|I_{1}\big\|_{L^{2}(\widetilde{M})}^{2}+\big\|I_{2}\big\|_{L^{2}(\widetilde{M})}^{2}
+\int_{\partial
\widetilde{M}}|\partial_{\nu}u|^{2}\,dS.
\end{equation}

\ms

{\bf Step 4}. In this step, we provide an
estimate for $\big\|\nabla
u\big\|_{L^{2}(\widetilde{M})}^{2}+\tau^{2}\big\||\nabla\varphi|u\big\|_{L^{2}(\widetilde{M})}^{2}$.

Let
\begin{equation}\label{4.11-2}
w=e^{\tau\varphi}(\partial_{x_{1}}-i\partial_{x_{2}})e^{-\tau\varphi}u.
\end{equation}
Then it follows that
\begin{equation}\label{9.20-eq4}
    \begin{aligned}
        \|w\|_{L^{2}(\widetilde{M})}^{2}=&\big\|u_{x_{1}}-\tau\varphi_{x_{1}}u\big\|_{L^{2}(\widetilde{M})}^{2}
        +\big\|u_{x_{2}}-\tau\varphi_{x_{2}}u\big\|_{L^{2}(\widetilde{M})}^{2}\\
        =&\big\|\nabla u\big\|_{L^{2}(\widetilde{M})}^{2}+\tau^{2}\big\||\nabla\varphi|u\big\|_{L^{2}(\widetilde{M})}^{2}
        -2\tau\int_{\widetilde{M}}(\nabla\varphi\cdot\nabla u)u\,dx.
    \end{aligned}
\end{equation}
Since $u\big|_{\partial\widetilde{M}}=0$ and
$\varphi$ is harmonic, we have
\begin{equation*}
    2\tau\int_{\widetilde{M}}(\nabla\varphi\cdot\nabla u)u\,dx
    =\tau\int_{\widetilde{M}}\mathrm{div}\,(u^{2}\nabla\varphi)\,dx
    -\tau\int_{\widetilde{M}}u^{2}\Delta\varphi\,dx=0.
\end{equation*}
This, together with \eqref{9.20-eq4}, implies
that
\begin{equation}\label{4.15}
    \big\|w\big\|_{L^{2}(\widetilde{M})}^{2}=\big\|\nabla u\big\|_{L^{2}(\widetilde{M})}^{2}+\tau^{2}\big\||\nabla\varphi|u\big\|_{L^{2}(\widetilde{M})}^{2}.
\end{equation}

From \eqref{4.11-2}, we see that
$$
\begin{cases} \ds
I_{1}= i(\partial_{x_{2}}+i\tau\varphi_{x_{1}})w=ie^{-i\tau\psi}\partial_{x_{2}}(e^{i\tau\psi}w),\\
\ns\ds I_{2}=
(\partial_{x_{1}}-i\tau\varphi_{x_{2}})w=e^{-i\tau\psi}\partial_{x_{1}}(e^{i\tau\psi}w).
\end{cases}
$$

Denote by
$$\Gamma_{-}=\big\{x\in\partial \widetilde{M}:\partial_{\nu}\varphi(x)<0\big\}.$$
Since $\varphi$ is harmonic, we have
$\int_{\partial
\widetilde{M}}\partial_{\nu}\varphi\,dS=0$.
Since $\varphi$ is a harmonic Morse function
with prescribed critical points
$\{p_{1},p_{2},\cdots,p_{m}\}$ in the interior
of $\widetilde{M}$, we know that $\varphi$ is
not a constant. Hence, $\Gamma_{-}$ is not
empty. Fix $\delta>0$ such that
$\Gamma_{\delta}=\{x\in\partial
\widetilde{M}:\partial_{\nu}\varphi(x)<-\delta\}$
is not empty.

Let $\psi$ be a conjugate function of $\varphi$,
i.e.,
$\Phi(x_{1},x_{2})=\varphi(x_{1},x_{2})+i\psi(x_{1},x_{2})$
is holomorphic in $\widetilde{M}$. By
\eqref{4.12}, we have for $\tau$ sufficiently
large, it holds that
\begin{equation}\label{4.14}
\begin{aligned}
\big\|w\big\|_{L^{2}(\widetilde{M})}^{2}&=\big\|e^{i\tau\psi}w\big\|_{L^{2}(\widetilde{M})}^{2}\\
&\leq
C\Big(\big\|\partial_{x_{1}}(e^{i\tau\psi}w)\big\|_{L^{2}(\widetilde{M})}^{2}
+
\big\|\partial_{x_{2}}(e^{i\tau\psi}w)\big\|_{L^{2}(\widetilde{M})}^{2}
+\int_{\Gamma_{\delta}}|\partial_{\nu}u|^{2}\,dS\Big)\\
&\leq
C\Big(\big\|ie^{-i\tau\psi}\partial_{x_{2}}(e^{i\tau\psi}w)\big\|_{L^{2}(\widetilde{M})}^{2}
+\big\|e^{-i\tau\psi}\partial_{x_{1}}(e^{i\tau\psi}w)\big\|_{L^{2}(\widetilde{M})}^{2}
-\tau\int_{\Gamma_{\delta}}\partial_{\nu}\varphi|\partial_{\nu}u|^{2}\,dS\Big)\\
&\leq
C\Big(\big\|I_{1}\big\|_{L^{2}(\widetilde{M})}^{2}+\big\|I_{2}\big\|_{L^{2}(\widetilde{M})}^{2}
-\tau\int_{\Gamma_{-}}\partial_{\nu}\varphi|\partial_{\nu}u|^{2}\,dS\Big).
\end{aligned}
\end{equation}

Combining \eqref{4.15} and \eqref{4.14}, we get
that
\begin{equation}\label{4.15-1}
\begin{aligned}
&\big\|\nabla u\big\|_{L^{2}(\widetilde{M})}^{2}+\tau^{2}\big\||\nabla\varphi|u\big\|_{L^{2}(\widetilde{M})}^{2}\\
&\leq
C\Big(\big\|I_{1}\big\|_{L^{2}(\widetilde{M})}^{2}+\big\|I_{2}\big\|_{L^{2}(\widetilde{M})}^{2}
-\tau\int_{\Gamma_{-}}\partial_{\nu}\varphi|\partial_{\nu}u|^{2}\,dS\Big).
\end{aligned}
\end{equation}

\ms

{\bf Step 5}. In this step, we prove that there
exists a constant $\tau_2>0$ such that for all
$\tau>\tau_2$, it holds that
\begin{equation}\label{4.16}
\tau\big\|u\big\|_{L^{2}(\widetilde{M})}^{2}\leq
C\big(\big\|u\big\|_{H^{1}(\widetilde{M})}^{2}+\tau^{2}\big\||\nabla\varphi|u\big\|_{L^{2}(\widetilde{M})}^{2}\big).
\end{equation}

In the following we use the notations
\begin{equation*}
\partial_{z}=\frac{1}{2}(\partial_{x_{1}}-i\partial_{x_{2}}),\quad
\partial_{\overline{z}}=\frac{1}{2}(\partial_{x_{1}}+i\partial_{x_{2}}).
\end{equation*}

For a critical point $q_{j}$ ($1\leq j\leq s$)
on $\partial\widetilde{M}$, we choose a
neighborhood $U_{j}\subset\widetilde{N}$ of
$q_{j}$ such that $q_{j}$ is the unique critical
point in $\overline{U_{j}}$, and
$\partial^{2}_{z}\varphi(x)\neq0$ for any $x\in
\overline{U_{j}}$. Here and in what follows, we
use the fact that the critical points of a Morse
function are non-degenerate. Since
$\partial\widetilde{M}$ is compact, we can
choose open sets $U_{s+1},\cdots,U_{s+r}$
satisfying that there are no critical points in
$\overline{U_{j}}$ ($s+1\leq j\leq s+r$), and
such that
$\partial\widetilde{M}\subset\bigcup_{j=1}^{s+r}U_{j}$.

For simplicity, we set $\mathcal{U}_{j}\deq
U_{j}\cap\widetilde{M}$, $\Gamma_{j,1}\deq
U_{j}\cap \partial\widetilde{M}$, and
$\Gamma_{j,2}\deq \partial U_{j}\cap
\widetilde{M}$ for $1\leq j\leq s+r$.

Next, we consider the critical points in the
interior of $\widetilde{M}$. For a critical
point $p_{k}$ ($1\leq k\leq m$), we choose a
neighborhood
$\mathcal{V}_{k}\subset\widetilde{M}$ of $p_{k}$
such that $p_{k}$ is the unique critical point
in it, and $\partial^{2}_{z}\varphi(x)\neq0$ for
any $x\in \overline{\mathcal{V}_{k}}$. We choose
also an open set
$\mathcal{V}_{0}\subset\widetilde{M}$ such that
there are no critical points in
$\overline{\mathcal{V}_{0}}$, and
$\widetilde{M}\subset
(\bigcup_{j=1}^{s+r}U_{j})\cup(\bigcup_{k=0}^{m}
\mathcal{V}_{k})$.

Let
$\{\zeta_{1},\cdots,\zeta_{s+r},\eta_{0},\eta_{1},\cdots,\eta_{m}\}$
be a smooth partition of unity on
$\widetilde{M}$, subordinate to the open sets
$\{U_{1},\cdots,U_{s+r},\mathcal{V}_{0},\mathcal{V}_{1},\cdots,\mathcal{V}_{m}\}$.
Then for $1\leq k\leq m$, since
$\mathrm{supp}\,\eta_{k}\subset
\mathcal{V}_{k}$, we can use integration by
parts to obtain
\begin{equation}\label{4.19-1}
\begin{aligned}
\tau\big\|\eta_{k}u\big\|_{L^{2}(\mathcal{V}_{k})}^{2}&=\tau\int_{\mathcal{V}_{k}}|\eta_{k}u|^{2}\,dx\leq C\tau\Big|\int_{\mathcal{V}_{k}}|\eta_{k}u|^{2}\overline{\partial_{z}^{2}\varphi}\,dx\Big|\\
&\leq C\tau\Big|\int_{\mathcal{V}_{k}}\partial_{\overline{z}}\Big[(\eta_{k}u)^{2}\Big]\overline{\partial_{z}\varphi}\,dx\Big| \\
&\leq C\big(\big\|\nabla(\eta_{k}u)\big\|_{L^{2}(\mathcal{V}_{k})}^{2}+\tau^{2}\big\||\nabla\varphi|\eta_{k}u\big\|_{L^{2}(\mathcal{V}_{k})}^{2}\big)\\
&\leq
C\big(\big\|u\nabla\eta_{k}\big\|_{L^{2}(\mathcal{V}_{k})}^{2}+
\big\|\eta_{k}\nabla
u\big\|_{L^{2}(\mathcal{V}_{k})}^{2}+\tau^{2}\big\||\nabla\varphi|\eta_{k}u\big\|_{L^{2}(\mathcal{V}_{k})}^{2}\big).
\end{aligned}
\end{equation}

For $1\leq j\leq s$, since
$\mathrm{supp}\,\zeta_{j}\subset U_{j}$, we can
use integration by parts and the boundary
condition $u\big|_{\Gamma_{j,2}}=0$ to obtain
\begin{eqnarray}\label{4.19-2}
\tau\big\|\zeta_{j}u\big\|_{L^{2}(\mathcal{U}_{j})}^{2} &=&\tau\int_{\mathcal{U}_{j}}|\zeta_{j}u|^{2}\,dx\leq C\tau\Big|\int_{\mathcal{U}_{j}}|\zeta_{j}u|^{2}\overline{\partial_{z}^{2}\varphi}\,dx\Big| \nonumber\\
&\leq&
C\tau\Big|\int_{\mathcal{U}_{j}}\partial_{\overline{z}}\Big[(\zeta_{j}u)^{2}\overline{\partial_{z}\varphi}\Big]\,dx
-\int_{\mathcal{U}_{j}}\partial_{\overline{z}}\Big[(\zeta_{j}u)^{2}\Big]\overline{\partial_{z}\varphi}\,dx\Big|\nonumber\\
&\leq&
C\tau\Big|\int_{\Gamma_{j,1}}(\zeta_{j}u)^{2}\overline{\partial_{z}\varphi}\,dS
+\int_{\Gamma_{j,2}}(\zeta_{j}u)^{2}\overline{\partial_{z}\varphi}\,dS
-\int_{\mathcal{U}_{j}}\partial_{\overline{z}}\Big[(\zeta_{j}u)^{2}\Big]\overline{\partial_{z}\varphi}\,dx\Big|\\
&\leq& C\tau\Big|\int_{\mathcal{U}_{j}}\partial_{\overline{z}}\Big[(\zeta_{j}u)^{2}\Big]\overline{\partial_{z}\varphi}\,dx\Big| \nonumber\\
&\leq& C\big(\big\|\nabla(\zeta_{j}u)\big\|_{L^{2}(\mathcal{U}_{j})}^{2}+\tau^{2}\big\||\nabla\varphi|\zeta_{j}u\big\|_{L^{2}(\mathcal{U}_{j})}^{2}\big) \nonumber\\
&\leq&
C\big(\big\|u\nabla\zeta_{j}\big\|_{L^{2}(\mathcal{U}_{j})}^{2}+
\big\|\zeta_{j}\nabla
u\big\|_{L^{2}(\mathcal{U}_{j})}^{2}+\tau^{2}\big\||\nabla\varphi|\zeta_{j}u\big\|_{L^{2}(\mathcal{U}_{j})}^{2}\big).\nonumber
\end{eqnarray}

Since $|\nabla\varphi|\neq0$ in
$\overline{U_{j}}$ ($s+1\leq j\leq s+r$) and
$\overline{\mathcal{V}_{0}}$, we know that there
exists $\tau_1>0$ such that for all
$\tau>\tau_1$, it holds that
\begin{equation}\label{4.19-3}
\tau\big\|\zeta_{j}u\big\|_{L^{2}(\mathcal{U}_{j})}^{2}
\leq
C\big(\big\|u\nabla\zeta_{j}\big\|_{L^{2}(\mathcal{U}_{j})}^{2}
+\big\|\zeta_{j}\nabla
u\big\|_{L^{2}(\mathcal{U}_{j})}^{2}
+\tau^{2}\big\||\nabla\varphi|\zeta_{j}u\big\|_{L^{2}(\mathcal{U}_{j})}^{2}\big),\q
s+1\leq j\leq s+r,
\end{equation}
and
\begin{equation}\label{4.19-4}
\tau\big\|\eta_{0}u\big\|_{L^{2}(\mathcal{V}_{0})}^{2}
\leq
C\big(\big\|u\nabla\eta_{0}\big\|_{L^{2}(\mathcal{V}_{0})}^{2}
+\big\|\eta_{0}\nabla
u\big\|_{L^{2}(\mathcal{V}_{0})}^{2}
+\tau^{2}\big\||\nabla\varphi|\eta_{0}u\big\|_{L^{2}(\mathcal{V}_{0})}^{2}\big).
\end{equation}

It follows from \eqref{4.19-1}--\eqref{4.19-4}
that
\begin{equation}\label{9.20-eq3}
\begin{aligned}
\tau\big\|u\big\|_{L^{2}(\widetilde{M})}^{2}&=\tau\Big\|\big(\sum_{j=1}^{s+r}\zeta_{j}+\sum_{k=0}^{m}\eta_{k}\big)u\Big\|_{L^{2}(\widetilde{M})}^{2}\\
&\leq\sum_{j=1}^{s+r}\tau\big\|\zeta_{j}u\big\|_{L^{2}(\mathcal{U}_{j})}^{2}
+\sum_{k=0}^{m}\tau\big\|\eta_{k}u\big\|_{L^{2}(\mathcal{V}_{k})}^{2}\\
&\leq
C\sum_{j=1}^{s+r}\big(\big\|u\nabla\zeta_{j}\big\|_{L^{2}(\mathcal{U}_{j})}^{2}+
\big\|\zeta_{j}\nabla
u\big\|_{L^{2}(\mathcal{U}_{j})}^{2}
+\tau^{2}\big\||\nabla\varphi|\zeta_{j}u\big\|_{L^{2}(\mathcal{U}_{j})}^{2}\big)\\
&\q+C\sum_{k=0}^{m}\big(\big\|u\nabla\eta_{k}\big\|_{L^{2}(\mathcal{V}_{k})}^{2}+
\big\|\eta_{k}\nabla
u\big\|_{L^{2}(\mathcal{V}_{k})}^{2}
+\tau^{2}\big\||\nabla\varphi|\eta_{k}u\big\|_{L^{2}(\mathcal{V}_{k})}^{2}\big)\\
&\leq
C\big(\big\|u\big\|_{L^{2}(\widetilde{M})}^{2}+
\big\|\nabla u\big\|_{L^{2}(\widetilde{M})}^{2}
+\tau^{2}\big\||\nabla\varphi|u\big\|_{L^{2}(\widetilde{M})}^{2}\big).
\end{aligned}
\end{equation}

By taking $\tau>\tau_2\deq \max\{C,\tau_1\}$,
where $C$ is the constant appearing in the right
hand side of \eqref{9.20-eq3}, we obtain
\eqref{4.16} from \eqref{9.20-eq3} immediately.

\ms

{\bf Step 6}. In this step, we complete the
proof. Combining \eqref{4.11}, \eqref{4.15-1}
and \eqref{4.16}, we get that
\begin{equation}\label{4.17}
\begin{aligned}
&\tau\big\|u\big\|_{L^{2}(\widetilde{M})}^{2}+\big\|u\big\|_{H^{1}(\widetilde{M})}^{2}
+\tau^{2}\big\||\nabla\varphi|u\big\|_{L^{2}(\widetilde{M})}^{2}+
\big\|\partial_{\nu}u\big\|_{L^{2}(\Gamma_{0})}^{2}\\
\leq & C\big(\big\|\nabla
u\big\|_{L^{2}(\widetilde{M})}^{2}
+\tau^{2}\big\||\nabla\varphi|u\big\|_{L^{2}(\widetilde{M})}^{2}+
\big\|\partial_{\nu}u\big\|_{L^{2}(\partial\widetilde{M})}^{2}\big)\\
\leq&C\Big(\big\|I_{1}\big\|_{L^{2}(\widetilde{M})}^{2}+\big\|I_{2}\big\|_{L^{2}(\widetilde{M})}^{2}
-\tau\int_{\Gamma_{-}}\partial_{\nu}\varphi|\partial_{\nu}u|^{2}\,dS+\|\partial_{\nu}u\|_{L^{2}(\partial \widetilde{M})}^{2}\Big)\\
\leq & C\Big(\big\|e^{\tau\varphi}\Delta
v\big\|_{L^{2}(\widetilde{M})}^{2}
+\tau\int_{\partial\widetilde{M}\setminus\Gamma_{-}}\partial_{\nu}\varphi|\partial_{\nu}u|^{2}\,dS\Big)\\
\leq &C\big(\big\|e^{\tau\varphi}\Delta
v\big\|_{L^{2}(\widetilde{M})}^{2}
+\tau\big\|\partial_{\nu}u\big\|_{L^{2}(\Gamma)}^{2}\big).
\end{aligned}
\end{equation}

Finally, for $u=e^{\tau\varphi}v$, since
$v\big|_{\partial \widetilde{M}}=0$, we have
$\partial_{\nu}u=e^{\tau\varphi}\partial_{\nu}v$
on $\partial\widetilde{M}$. From \eqref{4.17},
we  obtain that
\begin{equation}\label{4.18}
\begin{aligned}
&\tau\big\|e^{\tau\varphi}v\big\|_{L^{2}(\widetilde{M})}^{2}+\big\|e^{\tau\varphi}v\big\|_{H^{1}(\widetilde{M})}^{2}
+\tau^{2}\big\||\nabla\varphi|e^{\tau\varphi}v\big\|_{L^{2}(\widetilde{M})}^{2}+
\big\|e^{\tau\varphi}\partial_{\nu}v\big\|_{L^{2}(\Gamma_{0})}^{2}\\
\leq &C\big(\big\|e^{\tau\varphi}\Delta
v\big\|_{L^{2}(\widetilde{M})}^{2}
+\tau\big\|e^{\tau\varphi}\partial_{\nu}v\big\|_{L^{2}(\Gamma)}^{2}\big).
\end{aligned}
\end{equation}
For $q\in L^{\infty}(\widetilde{M})$, an
analogous inequality as \eqref{3.8-1} shows that
\begin{equation}\label{4.19}
\big\|e^{\tau\varphi}(-\Delta+q)
v\big\|_{L^{2}(\widetilde{M})}^{2}\geq\frac{1}{2}\big\|e^{\tau\varphi}\Delta
v\big\|_{L^{2}(\widetilde{M})}^{2}-\big\|qe^{\tau\varphi}v\big\|_{L^{2}(\widetilde{M})}^{2}.
\end{equation}
Combining \eqref{4.18} and \eqref{4.19}, by
taking $\tau_0 > \max\{\tau_2,
3C\|q\|_{L^{\infty}(\widetilde{M})}\}$, where
$C$ is the constant in the left hand side of
\eqref{4.18},  we see that for any
$\tau>\tau_0$, the inequality \eqref{1.7} holds.
This completes the proof Theorem \ref{thm1.2}.
\endpf

\section{Applications in the Calder\'{o}n problem}

As applications of the Carleman estimates
established in the previous sections, we give
some results on the Calder\'{o}n problem with
partial data.

Let $(M,g)$ be a $n$-dimensional  $C^{3}$-smooth
compact Riemannian manifold with the
$C^{2}$-smooth boundary. Set
$$
H_{\Delta_{g}}(M)\=\Big\{u\in
L^{2}(M):\Delta_{g} u\in L^{2}(M)\Big\},
$$
with the norm
$\|u\|_{H_{\Delta_{g}}(M)}=\|u\|_{L^{2}(M)}+\|\Delta_{g}
u\|_{L^{2}(M)}$.

By a similar argument as in \cite{BU02}, we know
that  there is a well defined bounded trace
operator from $H_{\Delta_{g}}(M)$ to
$H^{-\frac{1}{2}}(\partial M)$ and a normal
derivative operator from $H_{\Delta_{g}}(M)$ to
$H^{-\frac{3}{2}}(\partial M)$. Consequently,
the following set is well-defined:
$$
\mathcal{H}_{g}(\partial
M)=\big\{u\big|_{\partial M}:u\in
H_{\Delta_{g}}(M)\big\}\subset
H^{-\frac{1}{2}}(\partial M).
$$
Moreover, if $u\in H_{\Delta_{g}}(M)$ and
$u\big|_{\partial M}\in H^{\frac{3}{2}}(\partial
M)$, then $u\in H^{2}(M)$ and
$\partial_{\nu}u\big|_{\partial M}\in
H^{\frac{1}{2}}(M)$.

Let $q\in L^{\infty}(M)$. Assume that $0$ is not
a Dirichlet eigenvalue of $-\Delta_{g}+q$ on
$M$. Following \cite{A17},   for $f_{2}\in
\mathcal{H}_{g}(\partial M)$, the Dirichlet
problem
\begin{equation}\label{5.1}
\begin{cases}
(-\Delta_{g}+q)u=0 & \text{in}\,\, M,\\
u=f_{2} & \text{on}\,\, \partial M
\end{cases}
\end{equation}
has a unique solution $u\in H_{\Delta_{g}}(M)$.
The DN map
$\Lambda_{g,q}:\mathcal{H}_{g}(\partial
M)\rightarrow H^{-\frac{3}{2}}(\partial M)$ is
defined by
$$
\Lambda_{g,q}(f_{2})=\partial_{\nu}u\big|_{\partial
M},\q \forall f\in \mathcal{H}_{g}(\partial M).
$$
Let $\Gamma_{D}$ and $\Gamma_{N}$ be two open
subsets of $\partial M$. Define the partial
Cauchy data set as follows:
\begin{equation*}
\mathcal{C}_{g,q}^{\Gamma_{D},\Gamma_{N}}\=\Big\{(u\big|_{\Gamma_{D}},\partial_{\nu}u\big|_{\Gamma_{N}})
:(-\Delta_{g}+q)u=0\,\,\text{in}\,\,M,u\in
H_{\Delta_{g}}(M),\mathrm{supp}(u\big|_{\partial
M})\subset\Gamma_{D}\Big\}.
\end{equation*}

In the rest of this section, for the sake of
brevity, we use the  notation $e$ to denote the
Euclidean metric.

The Calder\'{o}n problem with partial data is to
determine $q$ from the knowledge of
$\mathcal{C}_{g,q}^{\Gamma_{D},\Gamma_{N}}$ for
given $\Gamma_{D}$ and $\Gamma_{N}$.

As we proved the Carleman estimates for
dimension $n\geq3$ and $n=2$ separately before,
below we also describe the partial data results
for dimension $n\geq3$ and $n=2$ separately.

By utilizing the Carleman estimates presented in
Section 1, one can proceed with the usual
approach to construct the CGO (Complex
Geometrical Optics) solution for the equation
\eqref{5.1} (e.g.,
\cite{BU02,DKSU09,IUY10,KS14,KSU07}).  By means
of the CGO solution, one can also follow  some
standard argument to obtain the following
results.
\begin{theorem}\label{prop5.1}
Let $(M,g)$ be an admissible manifold and assume
that $q_{1},q_{2}\in C(M)$. If $\partial
M_{\tan}$ has zero measure in $\partial M$, and
if $\mathcal{C}_{g,q_{1}}^{\partial
M_{-},\partial
M_{+}}=\mathcal{C}_{g,q_{2}}^{\partial
M_{-},\partial M_{+}}$, then $q_{1}=q_{2}$.
\end{theorem}

\begin{theorem}\label{prop5.2}
Let $\Omega\subset\mathbb{R}^{n}$ ($n\geq3$) be
a bounded open set with $C^{2}$-smooth boundary.
Let $x_{0}\in
\mathbb{R}^{n}\setminus\overline{\mathrm{ch}(\Omega)}$,
where $\mathrm{ch}(\Omega)$ is the convex hull
of $\Omega$. Denote by
\begin{equation*}
\begin{aligned}
F(x_{0})&=\{x\in\partial\Omega:(x-x_{0})\cdot\nu(x)\leq0\},\\
B(x_{0})&=\{x\in\partial\Omega:(x-x_{0})\cdot\nu(x)\geq0\}
\end{aligned}
\end{equation*}
the front and back face of $\Omega$. Let
$\Gamma_{D},\Gamma_{N}$ be two open subsets of
$\partial\Omega$ with
$F(x_{0})\subset\Gamma_{D}$ and
$B(x_{0})\subset\Gamma_{N}$. If $q_{1},q_{2}\in
L^{\infty}(\Omega)$, and if
$\mathcal{C}_{e,q_{1}}^{\Gamma_{D},\Gamma_{N}}=\mathcal{C}_{e,q_{2}}^{\Gamma_{D},\Gamma_{N}}$,
then $q_{1}=q_{2}$.
\end{theorem}

Theorem \ref{prop5.1} is established based on
the Carleman estimate \eqref{1.4} in Corollary
\ref{cor1.2}, while Theorem \ref{prop5.2} relies
on \eqref{1.5} in Corollary \ref{cor1.3}.

\begin{remark}
For $C^{\infty}$-smooth manifold  with a
$C^{\infty}$-smooth boundary, the same
uniqueness results as stated in Theorems
\ref{prop5.1} and  \ref{prop5.2} have been
proved in \cite{KS13} and \cite{KSU07}
respectively.
\end{remark}

For the case of $n=2$,
%the assumption of smoothness for the Riemann surface with  boundary is not essential as it can be conformally transformed to achieve smoothness. Additionally, for the set of partial Cauchy data, this only results in a conformal change of the potential.
based on the Carleman estimate \eqref{1.7} in
Theorem \ref{thm1.2}, the uniqueness result is
as stated below.
\begin{proposition}\cite[Theorem 1.1]{GT11}\label{prop5.3}
Let $(M,g)$ be a compact Riemann surface with
boundary. Let $\Gamma$ be an open subset of
$\partial M$. If $q_{1},q_{2}\in
C^{1,\alpha}(M)$ for some $\alpha>0$  and
$\mathcal{C}_{g,q_{1}}^{\Gamma,\Gamma}=\mathcal{C}_{g,q_{1}}^{\Gamma,\Gamma}$,
then $q_{1}=q_{2}$.
\end{proposition}
%
%\begin{proposition}\cite[Theorem 1.1]{IUY10}\label{prop5.4}
% Let $\Omega\subset\mathbb{R}^{2}$ be a bounded open set with $C^{2}$-smooth boundary. Let $\Gamma$ be an open subset of $\partial\Omega$. If $q_{1},q_{2}\in C^{2,\alpha}(\overline{\Omega})$ for some $\alpha>0$, and if $\mathcal{C}_{q_{1}}^{\Gamma,\Gamma}=\mathcal{C}_{q_{1}}^{\Gamma,\Gamma}$, then $q_{1}=q_{2}$.
%\end{proposition}


\begin{thebibliography}{10}

\bibitem{AFG17}
P. Angulo-Ardoy, D. Faraco and L. Guijarro,
\newblock {\em Sufficient conditions for the existence of limiting Carleman weights},
\newblock  Forum Math. Sigma {\bf 5} (2017), Paper No. e7, 25.

\bibitem{A17}
Y. Assylbekov,
\newblock {\em Reconstruction in the partial data Calder\'{o}n problem on admissible manifolds},
\newblock Inverse Probl. Imaging {\bf 11} (2017), no. 3, 455--476.

\bibitem{BU02}
A. L. Bukhgeim and G. Uhlmann,
\newblock {\em Recovering a potential from partial Cauchy data},
\newblock Comm. Partial Differential Equations {\bf 27} (2002), 653--668.

\bibitem{C39}
T. Carleman,
\newblock {\em Sur un probl\`eme d'unicit\'{e} pur les syst\`emes d'\'{e}quations
aux d\'{e}riv\'{e}es partielles \`a deux
variables ind\'{e}pendantes},
\newblock   Ark. Mat., Astr. Fys. 26 B. {\bf17} (1939), 1--9.

\bibitem{DKSU09}
D. Dos Santos Ferreira, C. E. Kenig, M. Salo and
G. Uhlmann,
\newblock {\em Limiting Carleman weights and anisotropic inverse problems},
\newblock  Invent. Math. {\bf 178} (2009), no. 1, 119--171.

%\bibitem{F91}
%O. Forster,
%\newblock {\em Lectures on Riemann surfaces},
%\newblock Graduate Texts in Mathematics {\bf 81}, Springer, New York, 1991.

\bibitem{FLZ19}
X. Fu, Q. L\"{u} and X. Zhang,
\newblock {\em Carleman estimates for second order partial differential
operators and applications, a unified approach},
\newblock   Springer, Cham, 2019.

\bibitem{GT11}
C. Guillarmou and L. Tzou,
\newblock {\em Calder\'{o}n inverse problem with partial data on {R}iemann surfaces},
\newblock   Duke Math. J. {\bf 158} (2011), no. 1, 83--120.

\bibitem{HT13}
B. Haberman and D. Tataru,
\newblock {\em Uniqueness in Calder\'{o}n's problem with Lipschitz conductivities},
\newblock Duke Math. J. {\bf 162} (2013), no. 3, 496--516.

\bibitem{H99}
E. Hebey,
\newblock {\em Nonlinear analysis on manifolds: {S}obolev spaces and
inequalities},
\newblock Courant Lecture Notes in Mathematics, vol. 5, New York University, Courant Institute of Mathematical Science, New York, American Mathematical Society, Providence, RI, 1999.

\bibitem{IUY10}
O. Y. Imanuvilov, G. Uhlmann and M. Yamamoto,
\newblock {\em The Calder\'{o}n problem with partial data in two dimensions},
\newblock J. Amer. Math. Soc. {\bf  23} (2010), no. 3, 655--691.

\bibitem{J17}
J. Jost,
\newblock {\em Riemannian geometry and geometric analysis},
\newblock seventh ed., Universitext, Springer, Cham, 2017.

%\bibitem{KS14}
%C. E. Kenig and M. Salo,
%\newblock {\em Recent progress in the Calder\'{o}n problem with partial data},
%\newblock Inverse Problems and Applications, 193--222, Contemp. Math. {\bf 615}, Amer. Math. Soc., Providence, RI, 2014.

\bibitem{KS13}
C. E. Kenig and M. Salo,
\newblock {\em The Calder\'{o}n problem with partial data on manifolds and applications},
\newblock Anal. PDE {\bf 6} (2013), no. 8, 2003--2048.

\bibitem{KS14} C. E. Kenig and M. Salo,
\newblock {\em Recent progress in the Calder\'{o}n problem with partial data}, in Inverse problems and applications, 193--222.
Contemp. Math., 615 American Mathematical
Society, Providence, RI, 2014

\bibitem{KSU07}
C. E. Kenig, J. Sj\"{o}strand and G. Uhlmann,
\newblock {\em The {C}alder\'{o}n problem with partial data},
\newblock Ann. of Math. (2) {\bf 165} (2007), no. 2, 567--591.

\bibitem{KkU16}
K. Krupchyk and G. Uhlmann,
\newblock {\em The {C}alder\'{o}n problem with partial data for conductivities with
$3/2$ derivatives},
\newblock Comm. Math. Phys. {\bf 348} (2016), no. 1, 185--219.

\bibitem{LRLR22}
J. Le Rousseau, G. Lebeau and L. Robbiano,
\newblock {\em Elliptic Carleman estimates and applications to stabilization and controllability. Vol. II. General boundary conditions on Riemannian manifolds},
\newblock   Progress in Nonlinear Differential Equations and their Applications, vol. 98, Birkh\"{a}user/Springer, Cham, 2022.

\bibitem{MT02}
R. Mazzeo and M. Taylor,
\newblock {\em Curvature and uniformization},
\newblock Israel J. Math. {\bf 130} (2002), 323--346.


\end{thebibliography}
\end{document}